\documentclass[12pt]{amsart}
\begin{document}
\numberwithin{equation}{section}

\def\1#1{\overline{#1}}
\def\2#1{\widetilde{#1}}
\def\3#1{\widehat{#1}}
\def\4#1{\mathbb{#1}}
\def\5#1{\frak{#1}}
\def\6#1{{\mathcal{#1}}}

\def\C{{\4C}}
\def\R{{\4R}}
\def\N{{\4N}}
\def\Z{{\4Z}}

\title[Obstructions to embeddability]{Obstructions to embeddability into hyperquadrics and explicit examples}
\author[D. Zaitsev]{Dmitri Zaitsev}
\address{D. Zaitsev: School of Mathematics, Trinity College Dublin, Dublin 2, Ireland}
\email{zaitsev@maths.tcd.ie}
\thanks{The author was supported in part by the RCBS grant of Trinity College Dublin 
and by the Science Foundation Ireland}

\maketitle
\begin{abstract} 
We give series of explicit examples of Levi-nondegenerate real-analytic hypersurfaces in complex spaces
that are not transversally holomorphically embeddable into hyperquadrics of any dimension.
For this, we construct invariants attached to a given hypersurface that serve as obstructions to embeddability.
We further study the embeddability problem for real-analytic submanifolds of higher codimension and answer a question
by Forstneri\v c.
\end{abstract} 
\tableofcontents

\def\Label#1{\label{#1}}


\def\cn{{\C^n}}
\def\cnn{{\C^{n'}}}
\def\ocn{\2{\C^n}}
\def\ocnn{\2{\C^{n'}}}


\def\dist{{\rm dist}}
\def\const{{\rm const}}
\def\rk{{\rm rank\,}}
\def\id{{\sf id}}
\def\aut{{\sf aut}}
\def\Aut{{\sf Aut}}
\def\CR{{\rm CR}}
\def\GL{{\sf GL}}
\def\Re{{\sf Re}\,}
\def\Im{{\sf Im}\,}
\def\wt{{\sf wt}\,}
\def\supp{{\rm supp}\,}
\def\span{{\rm span}\,}

\def\codim{{\rm codim}}
\def\crd{\dim_{{\rm CR}}}
\def\crc{{\rm codim_{CR}}}

\def\phi{\varphi}
\def\eps{\varepsilon}
\def\d{\partial}
\def\a{\alpha}
\def\b{\beta}
\def\g{\gamma}
\def\G{\Gamma}
\def\D{\Delta}
\def\Om{\Omega}
\def\k{\kappa}
\def\l{\lambda}
\def\L{\Lambda}
\def\z{{\bar z}}
\def\w{{\bar w}}
\def\Z{{\1Z}}
\def\t{\tau}
\def\th{\theta}

\emergencystretch15pt
\frenchspacing

\newtheorem{Thm}{Theorem}[section]
\newtheorem{Cor}[Thm]{Corollary}
\newtheorem{Pro}[Thm]{Proposition}
\newtheorem{Lem}[Thm]{Lemma}

\theoremstyle{definition}\newtheorem{Def}[Thm]{Definition}

\theoremstyle{remark}
\newtheorem{Rem}[Thm]{Remark}
\newtheorem{Exa}[Thm]{Example}
\newtheorem{Exs}[Thm]{Examples}

\def\bl{\begin{Lem}}
\def\el{\end{Lem}}
\def\bp{\begin{Pro}}
\def\ep{\end{Pro}}
\def\bt{\begin{Thm}}
\def\et{\end{Thm}}
\def\bc{\begin{Cor}}
\def\ec{\end{Cor}}
\def\bd{\begin{Def}}
\def\ed{\end{Def}}
\def\br{\begin{Rem}}
\def\er{\end{Rem}}
\def\be{\begin{Exa}}
\def\ee{\end{Exa}}
\def\bpf{\begin{proof}}
\def\epf{\end{proof}}
\def\ben{\begin{enumerate}}
\def\een{\end{enumerate}}
\def\beq{\begin{equation}}
\def\eeq{\end{equation}}

\section{Introduction}
The celebrated Whitney and Nash theorems about embeddings
of general smooth and Riemannian manifolds into their models (affine and euclidean spaces)
provide important tools for studying geometry of these manifolds.
The corresponding embeddability phenomena for Levi-nondegenerate real hypersurfaces in $\C^n$ (with models being the hyperquadrics in view of the Chern-Moser theory \cite{CM}) proves to be more difficult: On the one hand, Webster \cite{W78} showed
that any Levi-nondegenerate {\em real-algebraic} hypersurface 
is holomorphically (and transversally, see \S\ref{main-sect}) embeddable into a Levi-nondegenerate hyperquadric (of suitable dimension depending on the hypersurface).
On the other hand, it has been known since the work of Forstneri\v c \cite{Fo86} and Faran \cite{Fa88} that there exist strongly pseudoconvex real-analytic hypersurfaces not admitting any holomorphic embedding into a sphere in a complex space of any dimension. (More recently Forstneri\v c \cite{Fo04} extended these results to embeddings into hyperquadrics.)
In fact, such nonembeddable hypersurfaces 
have been shown to form a set of the first category in a suitable natural topology.
Despite of this, it appears to be difficult to obtain {\em explicit examples} of such hypersurfaces, none of which seems to be known.The given proofs used non-constructive Baire category arguments and did not lead to concrete examples.

It is one of the goals of the present paper to give explicit examples of nonembeddable real-analytic hypersurfaces of this kind.  Such hypersurfaces, in particular, cannot be algebraic nor even biholomorphically equivalent to algebraic ones in view of the mentioned result of Webster.
Thus we have to involve infinite power series:

\bt\Label{ex2}
Any hypersurface in $\C^2$ given by a convergent power series of the form
\begin{equation}
\Im w = z\bar z + \Re \sum_{k\ge 2} a_kT z^k \bar z^{(k+2)!}
+ \Re \sum b_{kml} z^k \bar z^m (\Re w)^l,
\end{equation}
where $a_k\ne 0$ for all $k$
and the second sum ranges over all $k,m,l$ satisfying 
$k,m\ge 2$,  $k\le (m+l)!$  and $m\le (k+l)!$,
 is not holomorphically embeddable into a sphere of any dimension.
More generally (see Remark~\ref{str} below), it is not transversally holomorphically embeddable into a hyperquadric of any dimension.
\et
As a special case we have, for instance, the following explicit nonembeddable example:
\bc\Label{ex1}
The hypersurface given in $\C^2$ by
\begin{equation}\Label{m1}
\Im w = z\bar z + \Re \sum_{k\ge 2} z^k \bar z^{(k+2)!}, \quad (z,w)\in \C^2, \quad |z|<\eps,
\end{equation}
for any $0<\eps\le 1$ is not transversally holomorphically embeddable into a hyperquadric of any dimension.
\ec
A brief glance at \eqref{m1} reveals the basic nature of this hypersurface:
in the expansion $\sum P_k(\bar z) z^k$ of the right-hand side,
the degree of the polynomial $P_k$ grows rapidly with respect to $k$.
These polynomials arise as partial derivatives in $z$ at $z=0$ of the complexified defining function of the hypersurface. For hypersurfaces of general form, we evaluate the partial derivatives of the defining function along the Segre varieties (see below for a definition). The above mentioned growth condition on the degrees of polynomials is then replaced by a rational dependence relation for those partial derivatives (which need not be polynomials in general). The latter dependence relations arise as obstructions to embeddability into a sphere
(or a hyperquadric) (see \S\ref{main-sect} and \S\ref{lower}).
These obstructions
can also be restated in terms of invariants attached to a given hypersurface
 (see \S\ref{invariants}) that can be directly calculated for given examples.
We further show that the latter invariants serve as lower bounds for the minimum possible dimension of the sphere (or hyperquadric),
where the given hypersurface can be holomorphically embeddable (resp.\ transversally holomorphically embeddable).
Those minimum dimensions (so-called ``CR complexity'') appear to be important (see e.g.\  the  recent work
of Baouendi-Ebenfelt-Huang~\cite{BEH}),
but to our knowledge, no bounds for them have been previously known.

The nonembeddability into a hyperquadric also implies the nonembeddability
into any Levi-nondegenerate {\em real-algebraic} hypersurface as a consequence of the mentioned result by Webster.
Note that examples of strongly pseudoconvex real-analytic hypersurfaces that are not 
{\em biholomorphically equivalent} (rather than embeddable) to any real-algebraic one are known due
to Huang-Ji-Yau~\cite{HJY}.
See also Baouendi-Ebenfelt-Rothschild~\cite{BERbull} for an earlier non-pseudoconvex example. Gausier-Merker~\cite{GM} gave obstructions to biholomorphic equivalence to real-algebraic submanifolds for generic real-analytic submanifolds that are of the tube form, have minimal possible dimension of their infinitesimal CR automorphism algebra and are minimal and finitely nondegenerate. In Appendix~\ref{sec-eq} we briefly illustrate how our methods can be utilized to obtain further such obstructions for general generic submanifolds without any restrictions, that can be used to obtain other examples.
However, it remains open whether these examples are transversally embeddable into
real-algebraic hypersurfaces of higher dimension.
Remarkably, embeddings into {\em infinite-dimensional} spheres or hyperquadrics always exist,
see Lempert~\cite{L82,L90} and D'Angelo~\cite{D84,D93}.

Forstneri\v c \cite{Fo04} furthermore obtains results
on the nonembeddability of ``most'' generic submanifolds of higher codimension
into {\em real-algebraic} generic submanifolds of possibly higher dimension but the {\em same codimension}.
However, his method does not apply to embeddings of CR manifolds of higher codimensions
into spheres or hyperquadrics (having codimension $1$) and hence he raises the question (Problem~2.3)
whether also here the set of all embeddable manifolds forms a set of the first category.
We answer this question affirmatively with our method in Theorem~\ref{nowhere}.

The paper is organized as follows.
In \S\ref{main-sect} we collect some basic material and notation,
state one of the main results about obstructions to embeddability
for hypersurfaces and give large series of explicit nonembeddable
examples based on these results. 
In Remark~\ref{rem-main} we illustrate the sharpness of Theorem~\ref{main}
by comparing the conclusion with the Chern-Moser theory in the case of biholomorphic equivalence.
\S\ref{algebra} contains the algebraic core of the paper,
where we establish transformation rules for jets of submanifolds
and defining functions that may be of independent interest.
The most crucial and nontrivial parts are the weight estimates.
In \S\ref{appl-embeddings} we apply the abstract results from the previous section
to embeddings of hypersurfaces into hyperquadrics.
This leads to obstructions to embeddability of hypersurfaces into hyperquadrics in terms
of their complex defining equations of the form $w=Q(z,\bar z,\bar w)$
(with $(z,w)$ not necessarily being the normal coordinates in the sense of \cite{BERbook}).
In \S\ref{poly-relations} we prove Theorem~\ref{main} in a slightly more general form.
In \S\ref{invariants} we attach to every germ of a real-analytic hypersurface
a series of invariant integers that can be used as lower estimates 
for the dimension of a hyperquadric, where the hypersurface is transversally embeddable.
In \S\ref{lower} we obtain lower order obstructions than those given by Theorem~\ref{main}
in the cases when the CR dimension of the given hypersurface is high.
This extends the phenomenon revealed by the Chern-Moser theory
that the obstructions for $M\subset\C^{n+1}$ from being equivalent to a hyperquadric
are of lower order for $n\ge 2$ than for $n=1$.
A new tool developed here is that of distinguished submanifolds of the Segre varieties
that are invariantly associated with real-analytic hypersurfaces.
Finally in \S\ref{high} we extend some of our results to embeddings of submanifolds
of higher codimension into hyperquadrics and give an answer to the question of Forstneri\v c mentioned above.

\section{Preliminaries and further results}\Label{main-sect}

Recall that a {\em real hyperquadric} in $\C^{n+1}$ is a real hypersurface given by
\begin{equation}\Label{quadric}
\Im w = h(z,\bar z)
\end{equation}
in some linear coordinates $(z,w)\in \C^n\times \C$, where $h(z,\bar z)$ is a hermitian form.
 By a {\em holomorphic embedding} of a real submanifold $M\subset \C^n$
into a real submanifold $M'\subset \C^{n'}$ we mean
a holomorphic embedding $H$ of a neighborhood of $M$ in $\C^n$ into $\C^{n'}$
with $H(M)\subset M'$. An embedding $H$ is said to be {\em transversal} to $M'$
if $H_* (T_x\C^n) + T_{H(x)} M' = T_{H(x)} \C^{n'}$ whenever $x\in M$.
We say that $M$ is {\em transversally holomorphically embeddable} into  $M'\subset \C^{n'}$ if there exists a transversal holomorphic embedding of $M$ into $M'$.
The transversality assumption is used to avoid trivial embeddings of $M$
into complex affine subspaces inside $M'$.
In case $M$ and $M'$ are hypersurfaces, transversality of an embedding of $M$ into $M'$ 
also guarantees that the Levi form of $M$ coincides with the restriction of that of $M'$.

\br\Label{str}
An embedding of a submanifold $M$ of {\em positive CR dimension}
into a {\em strictly pseudoconvex hypersurface} $M'$ is automatically transversal.
(Recall that the {\em CR dimension} of $M$ at $p\in M$ is the complex dimension of the complex tangent space $T^c_pM:=T_pM\cap iT_pM$.)
Indeed, for $p\in M$, write $L\colon T^c_pM\times T^c_pM\to (T_pM/T^c_pM)\otimes \C$ 
for the Levi form and use the corresponding notation for $M'$.
If $H$ is any holomorphic map of a neighborhood of $M$ in $\C^n$ into $\C^{n'}$ with $H(M)\subset M'$,
then one has $H_*(T^c_pM)\subset T^c_{H(p)}M'$ and $H_* L(u,u)=L'(H_*u,H_*u)$ for $u\in T^c_pM$.
If $H$ is an embedding, one has $H_*u\ne 0$ for $u\ne 0$.
Then, since $M'$ is strongly pseudoconvex, one has $L'(H_*u,H_*u)\ne0$ and therefore $H_* L(u,u)\ne 0$.
The latter is a condition in $(T_{H(p)} M'/T^c_{H(p)}M') \otimes \C$,
which easily implies the transversality.

The reader is referred to Ebenfelt-Rothschild \cite{ER} for a more general and detailed analysis of transversality of holomorphic maps between CR manifolds of the same dimension and to Baouendi-Ebenfelt-Rothschild \cite{BER07} for holomorphic maps between hypersurfaces of different dimensions.
\er

Let $M\subset \C^{n+1}$ be a real-analytic hypersurface with a fixed reference point $p\in M$
that we shall assume to be $0$ for simplicity.
We choose local holomorphic coordinates $Z=(z,w)\in \C^n\times \C$ 
defined in a neighborhood of $0$ and vanishing at $0$, 
together with a real-analytic defining function $\rho(Z,\bar Z)$ of $M$
(that we think of as a  power series in $(Z,\bar Z)$ convergent in a neighborhood of the origin)
satisfying $\rho_w(0,0)\ne 0$. Recall that the {\em Segre variety}
of a point $\zeta\in\C^{n+1}$, associated to $M$,
is defined by $S_{\zeta} := \{Z : \rho (Z,\bar\zeta)=0\}$
(see e.g.\ \cite{W77} for basic facts about Segre varieties).
In particular, we shall consider the Segre variety $S_0$ corresponding to the reference point $0$.
{\em We shall always assume that $S_0$ is connected}.
We use the subscript notation (sometimes separated by commas) for the partial derivatives,
i.e.\ 
$$\rho_{z^\b w^s} =\rho_{z^\b,w^s} 
= \frac{\d^{|\b|+s} \rho}{\d z_1^{\b_1} \cdots \d z_n^{\b_n}\d w^s},$$
for a multiindex $\b = (\b_1,\ldots,\b_n)\in\N^n$ and an integer $s\in\N$
(where we keep the usual convention that $\N$ begins with $0$).
To the derivatives of $\rho$, we assign the following weights 
that will appear natural in the transformation rules below:
\begin{equation}\Label{wt40}
\wt \rho_{z^\a w^s}:= 2|\a|+s-1, \quad \a\in \N^n, \; s\in \N.
\end{equation}
Thus the derivative $\rho_w$ is the only one with weight $0$
and hence not contributing to the total weight.
Therefore we shall count it separately, denoting by $\deg_{\rho_w} P$
the degree of the polynomial $P$ in the variable $\rho_w$
(i.e.\ the maximal power of $\rho_w$ that appears in $P$).

\bt\Label{main}
Let $M\subset \C^{n+1}$ be a real-analytic hypersurface through $0$ given by $\rho(Z,\bar Z)=0$
with $\rho_w(0,0)\ne 0$.
Suppose that $M$ is transversally holomorphically embeddable 
into a hyperquadric in $\C^{n+m+1}$.
Then for any set of $m+1$ multiindices $\a_j\in \N^n$, $|\a_j|\ge 2$, $j=1,\ldots,m+1$,
there exists an integer $k$ with $K:=\{j: |\a_j| = k\}\ne \emptyset$
such that the partial derivatives of $\rho$ satisfy a relation of the form
\begin{equation}\Label{relate0}
\sum_{j\in K} P_j(\rho_{z^\b w^s}(0,\bar\zeta))\, \rho_{z^{\a_j}}(0,\bar\zeta) = 
R(\rho_{z^\b w^s}(0,\bar\zeta)), 
\quad\zeta\in S_0,
\end{equation}
where $P_j(\rho_{z^\b w^s}(0,\bar\zeta))$ and $R(\rho_{z^\b w^s}(0,\bar\zeta))$ are polynomials
in the partial derivatives $\rho_{z^\b w^s}(0,\bar\zeta)$ with $0<|\b|+s\le k$, $0<|\b|<k$,
such that not all $P_j(\rho_{z^\b w^s}(0,\bar\zeta))$ identically vanish in $\zeta\in S_0$.
Moreover, $P_j$ and $R$ can be chosen satisfying in addition the following weight and degree estimates:
\begin{equation}\Label{wt-est0}
\begin{array}{lll}
&\wt P_j\le (2k-2)(l-1),\quad &\wt R\le (2k-2)l+1, \\
&\deg_{\rho_w} P_j \le (2k-2)l, \quad &\deg_{\rho_w} R \le (2k-2)l+1,
\end{array}
\end{equation}
where $l\ge1$ is the number of all $j$ with $|\a_j|\le k$.
\et

In fact, we shall prove a more general version of Theorem~\ref{main} in the form Theorem~\ref{main1} below,
where the Segre variety $S_0$ is replaced by any irreducible subvariety through $0$.
In the special case when $M$ is {\em rigid} (in the sense of Baouendi-Rothschild), Theorem~\ref{main}
can be stated in the following simplier form without referring to Segre varieties. The proof is straightforward.

\bc\Label{main-rigid}
Let $M\subset \C^{n+1}$ be a real-analytic hypersurface through $0$ given
in its rigid form by $\Im w=\phi(z,\bar z)$, $(z,w)\in\C^n\times\C$.
Suppose that $M$ is transversally holomorphically embeddable 
into a hyperquadric in $\C^{n+m+1}$.
Then for any set of $m+1$ multiindices $\a_j\in \N^n$, $|\a_j|\ge 2$, $j=1,\ldots,m+1$,
there exists an integer $k$ with $K:=\{j: |\a_j| = k\}\ne \emptyset$
such that the partial derivatives of $\phi$ satisfy a relation of the form
\begin{equation}\Label{relate10'}
\sum_{j\in K} P_j(\phi_{z^\b}(0,\bar\chi))\, \phi_{z^{\a_j}}(0,\bar\chi) = 
R(\phi_{z^\b}(0,\bar\chi)), 
\quad\chi\in \C^n,
\end{equation}
where $P_j(\phi_{z^\b}(0,\bar\chi))$ and $R(\phi_{z^\b}(0,\bar\chi))$ are polynomials
in the partial derivatives $\rho_{z^\b}(0,\bar\chi)$ with $0<|\b|<k$,
such that not all $P_j(\rho_{z^\b}(0,\bar\chi))$ identically vanish in $\chi$.
Moreover, $P_j$ and $R$ can be chosen satisfying in addition the following weight estimates:
\begin{equation}\Label{wt-est0'}
\wt P_j\le (2k-2)(l-1),\quad \wt R\le (2k-2)l+1,
\end{equation}
where $l\ge1$ is the number of all $j$ with $|\a_j|\le k$.
\ec

\br\Label{rem-main}
We here consider the special case $m=0$, where the conclusion of Theorem~\ref{main}
can be compared with that of the Chern-Moser theory \cite{CM}.
For $m=0$, Theorem~\ref{main} gives obstructions preventing $M$
from being (locally) biholomorphically  {\em equivalent} to a hyperquadric.
Of course, the full set of such obstructions is known due to 
the Chern-Moser normal form \cite{CM},
whose actual computation, however, may be hard in concrete cases.
On the other hand, Theorem~\ref{main} may be applied directly in given coordinates
instead of the normal coordinates obtained through the Chern-Moser normalization.
For instance, for a single multiindex $|\a|=2$, Theorem~\ref{main} yields
(with $k=2$, $l=1$) a relation
\begin{equation}\Label{m0}
\rho_{z^\a} = R(\rho_{z^\b}, \rho_{z^\b w},
\rho_{w^2},\rho_w), \quad |\b|=1,
\end{equation}
where $R$ is a polynomial of weight $\le 3$
and all derivatives are evaluated at $(0,\bar\zeta)$, $\zeta\in S_0$.
Thus, if \eqref{m0} is not satisfied, $M$ is not equivalent to a hyperquadric.
In particular, if $M$ is in its Chern-Moser normal form \cite{CM}, we have
\begin{equation*}
\rho_{z^\b w}(0,\bar\zeta)\equiv \rho_{w^2}(0,\bar\zeta) \equiv 0,
\quad \rho_w(0,\bar\zeta)\equiv\const,
\quad \rho_{z^\b}(0,\bar\zeta) \text{ is linear in }\bar \zeta.
\end{equation*}
Then \eqref{m0} means that any $2$nd order derivative $\rho_{z^\a}(0,\bar\zeta)$
is a polynomial in $\bar\zeta$ of degree $\le 3$.
We now compare this with the normal form $M=\{\Im w = \sum a_{\a\mu s} z^\a \bar z^\mu (\Re w)^s\}$,
where normalization conditions are imposed, in particular, 
on the coefficients $a_{\a\mu 0}$ with $|\mu|\le 3$.
Here \eqref{m0} means the vanishing of the coefficients $a_{\a\mu 0}$ with $|\mu|\ge 4$,
which are exactly the free coefficients that appear in the normal form and
hence have to vanish in order for $M$ to be equivalent to a hyperquadric. 
Thus the estimates given by \eqref{wt-est0} are sharp in this case.
\er

Based on Theorem~\ref{main}, one can obtain explicit examples
of hypersurfaces that are not transversally embeddable into hyperquadrics of certain dimensions
or into hyperquadrics of any dimension.

\begin{proof}[Proof of Theorem \ref{ex2}]
We write 
$$\rho(Z,\bar Z):= -\Im w + z\bar z + \Re \sum_{k\ge 2} z^k \bar z^{(k+2)!}+ \Re \sum b_{kml} z^k \bar z^m (\Re w)^l
$$
with the second sum ranging as in the assumption.
Then 
 $M$ is given by $\rho(Z,\bar Z)=0$, we have $S_0=\{w=0\}$ for the Segre variety of $0$ and
$\rho_w=-\frac1{2i}$, $\rho_z(0,\bar \zeta) = \bar \chi$, where $\zeta=(\chi,0)\in S_0\subset \C^n\times\C$.
Furthermore, $\rho_{z^a}(0,\bar \zeta)$ is a polynomial in $\bar \chi$ of degree $(a+2)!$ for every $a\ge2$,
and every other derivative $\rho_{z^a w^b}(0,\bar \zeta)$, $b\ge1$, is a polynomial in $\bar\chi$
of degree $\le (a+b)!$.
By contradiction, assume that $M$ is transversally holomorphically embeddable into a hyperquadric in some $\C^{2+m}$.
Then, in view of Theorem~\ref{main} applied to $\a_j=j+1$, $j=1,\ldots,m+1$, 
there is a $k\ge 2$ with $K:=\{j: |\a_j| = k\}=\{k-1\}$ and a 
relation \eqref{relate0} with $P_{k-1}(\rho_{z^\b w^s}(0,\bar\zeta))$ and 
$R(\rho_{z^\b w^s}(0,\bar\zeta))$ satisfying \eqref{wt-est0}.
In particular, we have $\wt R \le (2k-2)(k-1)+1$ in view of $l=k-1$.
Since $\wt \rho_{z^aw^b}=2a+b-1$,
we have 
$\deg\rho_{z^a w^b}(0,\bar\zeta)\le \frac{(k+1)!}{2k-3}\wt \rho_{z^a w^b}$
for every $a,b$ satisfying $a+b\le k$, $a<k$. Then
it  follows that $R(\rho_{z^\b w^s}(0,\bar\zeta))$
is a polynomial in $\bar\zeta$ whose degree does not exceed
$$\frac{(k+1)!}{2k-3} \wt R \le \frac{(k+1)!}{2k-3}  ((2k-2)(k-1)+1) < (k+2)!.$$
This is a contradiction with \eqref{relate0}
since $\rho_{z^k}(0,\bar\zeta))$ is of degree precisely $(k+2)!$.
The proof is complete.
\end{proof}

\section{Some algebraic operations with multilinear functions and transformation formulas}
\Label{algebra}

\subsection{An algebra of symmetric multilinear functions}
We fix a finite-dimensional complex vector space $V$
and denote by $\6P_d$, $d=0,1,\ldots$, the space of all symmetric $d$-linear functions
$$p\colon V\times \cdots \times V = V^d\to \C,$$
and by $\6P:=\oplus_d \6P_d$ the corresponding graded direct sum.
In case $d=0$ we set $\6P_0:=\C$, i.e.\ ``$0$-linear'' functions are identified with complex numbers.
We write $\deg p = d$ for $p\in \6P_d\setminus\{0\}$.
There is a standard one-to-one correspondence between $\6P_d\setminus\{0\}$
and the homogeneous polynomials on $V$ of degree $d$ obtained by associating to every $p\in\6P_d$
its evaluation $p(v,\ldots,v)$.
Then the product of polynomials induces a natural product on $\6P$.
However, it will be more convenient for our purposes to consider
{\em another product} on $\6P$ that differs from the mentioned ``polynomial product'' 
by certain additional factors depending on the degree.
As a result, there will be less additional factors in the transformation formulas below.

The product we consider here can be defined as follows.
For $p_1(v_1,\ldots,v_{d_1})\in \6P_{d_1}$ and $p_2(v_1,\ldots,v_{d_2})\in \6P_{d_2}$ 
define
\begin{equation}\Label{multiplication}
(p_1\cdot p_2)(v_1,\ldots,v_{d_1+d_2}):= 
\sum p_1(v_{i_1},\ldots,v_{i_{d_1}})\, p_2(v_{j_1},\ldots,v_{j_{d_2}}),
\end{equation}
where the summation is taken over all possible (disjoint) partitions 
$$\{1,\ldots,d_1+d_2\} = \{i_1,\ldots,i_{d_1}\} \cup \{j_1,\ldots,j_{d_2}\}.$$
It is easy to see that $p_1\cdot p_2$ so defined is again symmetric in its arguments
and hence $p_1\cdot p_2\in \6P_{d_1+d_2}$.
It is furthermore easy to check that this operation of multiplication together with
the usual addition makes $\6P$ a commutative associative graded $\C$-algebra with unit $1\in \6P_0$.

\be
For $V=\C$, $p_1(x_1)=x_1\in\6P_1$ and $p_2(x_1,x_2)=x_1x_2\in\6P_2$, we have 
$$(p_1\cdot p_2)(x_1,x_2,x_3) = p_1(x_1) p_2(x_2,x_3) + p_1(x_2) p_2(x_3,x_1) + p_1(x_3) p_2(x_1,x_2)\in\6P_3,$$
whereas the ``polynomial product'' would give $\frac13 (p_1\cdot p_2)$.
\ee

We next consider an operation of substitution (or composition).
Let $A_j\colon  V^{\nu_j}\to V$, $\nu_j\ge0$, $j=1,\ldots,m$,
be a set of maps, where each $A_j$ is symmetric $\nu_j$-linear.
As before, a ``$0$-linear'' map $A_j\colon V^0\to V$ means by definition a vector in $V$.
We shall write $(A_1,\ldots,A_m)=A_{\nu_1,\ldots,\nu_m}$ indicating
the degrees as subscripts.
For $p\in \6P_d$ with $d\ge m$, we then define the ``substitution'' 
$p\circ A_{\nu_1,\ldots,\nu_m}\in \6P_{d-m+\nu_1+\cdots+\nu_m}$ as follows:
\begin{multline}\Label{}
(p\circ A_{\nu_1,\ldots,\nu_m})(v_1,\ldots,v_{d-m+\nu_1+\cdots+\nu_m}):= \\
\sum p(A_1(v_{a^1_1},\ldots,v_{a^1_{\nu_1}}), \ldots, A_m(v_{a^m_1},\ldots,v_{a^m_{\nu_m}}),
v_{b_1},\ldots,v_{b_{d-m}}),
\end{multline}
where the summation is taken over all possible partitions
$$\{1,\ldots,d-m+\nu_1+\cdots+\nu_m\}= 
\{a^1_1,\ldots,a^1_{\nu_1}\} \cup \cdots \cup \{a^m_1,\ldots,a^m_{\nu_m}\}
\cup \{b_1,\ldots,b_{d-m}\}.$$
Again it is easy to see that the result is symmetric in its arguments and hence is in 
$\6P_{d-m+\nu_1+\cdots+\nu_m}$.
It will also be convenient to allow the case $m=0$, i.e.\ consider the substitution 
of the empty set $\emptyset$ of maps $A_j$ into $p$, 
where we define $p\circ \emptyset := p$.

What is the result of the substitution operation applied twice?
It is not difficult to see that such repeated substitution
is actually a sum of single substitutions.
More precisely, we have the following elementary lemma, the proof of which is straightforward.
We use the notation $\deg A_j = \nu_j$ if $A_j\colon  V^{\nu_j}\to V$ is $\nu_j$-linear.

\bl\Label{composition}
Let $p':=p\circ A_{\nu_1,\ldots,\nu_m}$ be as before and let 
$B_s\colon  V^{\mu_s}\to V$, $s=1,\ldots,l$, be another collection of symmetric multilinear maps
such that the composition $p'':=p'\circ B_{\mu_1,\ldots,\mu_l}$ is defined (i.e.\ $l\le \deg p'$).
Then $p''$ is a finite sum of terms of the form  $p\circ C_{\lambda_1,\ldots,\lambda_s}$,
each with suitable multilinear maps $C_j\colon V^{\lambda_j}\to V$, satisfying
$$\lambda_1+\cdots +\lambda_r \le (\nu_1+\cdots +\nu_m) + (\mu_1+\cdots +\mu_l).$$
\el

\subsection{Transformation of submanifolds jets via embeddings}\Label{j-trans}
Our goal here is to obtain a relation formula between jets of complex submanifolds
and of their embeddings with explicit degree and weight estimates.
We consider a holomorphic embedding $H$
from a neighborhood of $0$ in $\C^{n+1}$ into $\C^{n+m+1}$
and split the coordinates as follows: $(z,w)\in \C^n\times \C$ and $(z',w')\in \C^n\times \C^{m+1}$.
Consider complex hypersurfaces $S$ in $\C^{n+1}$ passing through $0$
and their images $S'=H(S)\subset \C^{n+1}$, both represented as graphs of holomorphic functions
$w=Q(z)$ and $w'=Q'(z')$ respectively.
Thus we have the relation
\begin{equation}\Label{basic}
G(z,Q(z))=Q'(F(z,Q(z))).
\end{equation}
We want to express the derivatives of $Q'$ in terms of the derivatives of $Q$, $F$ and $G$.
In general, these expressions are rational but we shall make a first order assumption on $H$
making the relations polynomial.
Writing 
\begin{equation}\Label{split}
H(z,w)=(F(z,w),G(z,w)) \in \C^n\times \C^{m+1}
\end{equation} 
with respect to the chosen coordinates,
our main assumption is
\begin{equation}\Label{normalization}
F_z(0)=\id, \quad F_w(0)=0,
\end{equation}
where $\id$ stands for the identity $n\times n$ matrix.
We write $Q_{z^k}$ for the full $k$th derivative at $0$,
i.e.\ $Q_{z^k}$ is a $k$-linear function $(\C^n)^k\to \C$ given in terms of the partial derivatives by
\begin{equation}\Label{}
Q_{z^k}(v^1,\ldots,v^k):= \sum Q_{z_{j_1},\ldots,z_{j_k}}(0) v^1_{j_1} \ldots v^k_{j_k},
\end{equation}
where $v^s=(v^s_1,\ldots,v^s_n)\in\C^n$ and the summation is taken 
over all multiindices $(j_1,\ldots,j_k)\in \{1,\ldots,n\}^k$.
In case $k=0$ we set $Q_{z^0}:=1\in \C=\6P_0$.
Similar notation will be used for $G$:
\begin{equation}\Label{g-notation}
G_{z^k w^l}(v^1,\ldots,v^k):= \sum G_{z_{j_1},\ldots,z_{j_k},w^l}(0) v^1_{j_1} \ldots v^k_{j_k},
\end{equation}
where the full derivative is only taken with respect to $z$.
The derivatives of $F$ will be regarded in the same way but will be suppressed
in our transformation formula below, whereas the derivatives of $G$ will appear more explicitly.

We next introduce weights of the derivative terms as follows. We first set
\begin{equation}\Label{wt1}
\wt G_{z^s w^l}:= 2s+l-1, \quad \wt Q_{z^s}:= 2s-1,
\end{equation}
and then extend them to compositions by
\begin{equation}\Label{wt2}
\wt G_{z^s w^l}\circ A_{\nu_1,\ldots,\nu_a} := \wt G_{z^s w^l} + \nu_1+\ldots+\nu_a,
\quad \wt Q_{z^s}\circ A_{\nu_1,\ldots,\nu_a} := \wt Q_{z^s} + \nu_1+\ldots+\nu_a.
\end{equation}
That is, for every composition, the sum of the total degrees of the multilinear maps 
$A_{\nu_1},\ldots,A_{\nu_a}$ is simply added to the weight of $G_{z^s w^l}$ or $Q_{z^s}$.

\bp\Label{transformation}
Under the normalization assumption \eqref{normalization},
the full higher order derivatives of $Q$ and 
$Q'$ at $0$ are related by the formula
\begin{equation}\Label{transfQQ'}
Q'_{z'^k} = \sum (G_{z^s w^l}\circ A_{\nu_1,\ldots,\nu_a})\cdot 
(Q_{z^{s^1}}\circ B_{\nu_1^1,\ldots,\nu_{a^1}^1})\cdot\ldots\cdot
(Q_{z^{s^r}}\circ B_{\nu_1^r,\ldots,\nu_{a^r}^r}),
\end{equation}
where the summation is taken over all (finitely many) indices 
$s,l$, sets of indices $\{s^1,\ldots,s^r\}$ with $r\ge l$, 
and finitely many sets of multilinear maps 
$A_{\nu_1,\ldots,\nu_a}$ and $B_{\nu_1^j,\ldots,\nu_{a^j}^j}$ 
(including some of them or all being empty sets) depending only on $F$,
such that the degree of each term on the right-hand side of \eqref{transfQQ'}
equals $k=\deg Q'_{z'^k}$ and its weight does not exceed $2k-1 = \wt Q'_{z'^k}$.
Moreover, each term with the empty set of multilinear maps appears precisely once.
\ep

\bpf
We proceed by induction on $k$.
The case $k=1$ is easy and obtained by direct differentiating \eqref{basic} in $z$
and using the normalization \eqref{normalization}:
$$Q'_{z'}= G_z + G_w \cdot Q_z.$$

We now assume that \eqref{transfQQ'} holds for all $k<k_0$
and take the full $k_0$th derivatives of both sides in \eqref{basic} evaluated at $0$.
On the left-hand side we obtain the terms 
\begin{equation}\Label{g-term}
G_{z^s w^l}\cdot Q_{z^{s_1}}\cdot \ldots\cdot Q_{z^{s_l}},
\quad s+ s_1+\ldots + s_l = k_0,
\end{equation}
(with the number $l$ of factors $Q_{z^{s_r}}$ being equal to the
$w$-order in $G_{z^s w^l}$).
According to our definition of multiplication of multilinear maps \eqref{multiplication},
we obtain precisely one term of the form \eqref{g-term}
for each choice of $s,l$ and of a (possibly empty) set of indices $\{s_1,\ldots,s_l\}$. 
The weight of \eqref{g-term} is 
$$2s+l-1 + \sum_{r=1}^l (2s_r-1) = 2k_0-1,$$
as desired.
Similarly, on the right-hand side, the terms will be of the form
\begin{equation}\Label{f-terms}
Q'_{z'^k} \left( 
(F_{z^{s^1} w^{a^1}}\cdot Q_{z^{s^1_1}}\cdot \ldots\cdot Q_{z^{s^1_{a^1}}}),
\ldots,
(F_{z^{s^k} w^{a^k}}\cdot Q_{z^{s^k_1}}\cdot \ldots\cdot Q_{z^{s^k_{a^k}}})
\right),
\end{equation}
with $\sum_j (s^j+ s^j_1+\ldots + s^j_{a^j}) = k_0$, where we regard
$Q'_{z'^k}$ as before as a multilinear function with $k$ arguments.
Here $k\le k_0$ and there is precisely one term with $k=k_0$,
namely $Q'_{z'^{k_0}}$ itself, where we continue using our normalization \eqref{normalization}.
Thus we express $Q'_{z'^{k_0}}$ as the left-hand side minus the terms on the right-hand
side with $k<k_0$. For the latter terms we can use our induction hypothesis
that each $Q'_{z'^k}$, $k<k_0$, is already given by the formula \eqref{transfQQ'}.
Substituting it into \eqref{f-terms} and using Lemma~\ref{composition}
and the fact that each $Q_{z^{s^j_r}}$ is scalar,
we conclude that each term in \eqref{f-terms} with $k<k_0$ is
expressible as a finite sum of the terms in \eqref{transfQQ'}.
Clearly the total degree of each such term is always $k_0$.

It remains to show that the weight of each term does not exceed $2k_0-1$.
Each term in \eqref{f-terms} with $k<k_0$ 
arises as a composition of a term $p\in\6P_k^{m+1}$ in \eqref{transfQQ'}
of weight $\le 2k-1$ with $k$ multilinear maps
\begin{equation}\Label{maps}
(F_{z^{s^1} w^{a^1}}\cdot Q_{z^{s^1_1}}\cdot \ldots\cdot Q_{z^{s^1_{a^1}}}),
\ldots,
(F_{z^{s^k} w^{a^k}}\cdot Q_{z^{s^k_1}}\cdot \ldots\cdot Q_{z^{s^k_{a^k}}}).
\end{equation}
We first look at the extreme cases, where all $F$-derivatives in \eqref{maps} are $F_z=\id$
except one, which is either $F_{zw}\cdot Q_{z^{k_0-k}}$ or 
$F_{w^2}\cdot Q_z \cdot Q_{z^{k_0-k}}$.
The corresponding compositions are $(p\circ F_{zw}) \cdot Q_{z^{k_0-k}}$
and $(p\circ F_{w^2}) \cdot Q_z \cdot  Q_{z^{k_0-k}}$,
both having weight $\le (2k-1)+1 + 2(k_0-k)-1 = 2k_0-1$ in view of Lemma~\ref{composition}.
Note that by \eqref{normalization}, there is no term with $F_w$.

Our strategy to estimate the weights of general terms is to compare them with these extreme cases.
More precisely, we shall consider simple moves to pass from one term to another.
Our first move consists of raising the $z$-order $s$ in $F_{z^s w^r}$ by an integer $t$.
In order to keep the total degree constant, we decrease
by the same integer $t$ the order $l$ in some factor $Q_{z^l}$.
Since the increase contributes with $+t$ to the total weight,
whereas the decrease with $-2t$ in view of our rules \eqref{wt1}-\eqref{wt2},
we can only decrease the total weight that way.
Our second move raises the $w$-order in $F_{z^s w^r}$ by $r'$
and adds $r'$ new factors $Q_z$. Again, to keep the total degree constant,
we have to lower by $r'$ the order of $Q_{z^l}$.
Then the total weight increases by $r'$ and decreases by $2r'$, hence decreases in total.
Using these two moves we shall obtain any term with all maps in \eqref{maps} being $F_z=\id$
except one, being
\begin{equation}\Label{f-prod}
F_{z^s w^r}\cdot Q_{z^l}\cdot Q_z \ldots \cdot Q_z
\end{equation}
with appropriate integers and appropriate number of the first order factors $Q_z$.
Our next move exchanges derivative orders between the $Q$-factors here.
That is, keeping the total degree constant, we can decrease the order of a factor in \eqref{f-prod}
by an integer and simultaneously increase the order of another $Q$-factor by the same integer.
Clearly this move does not change the weight and allows us to obtain
any other term still having all but one maps in \eqref{maps} equal $F_z=\id$.

Our two last moves will exchange indices between different parentheses in \eqref{maps}.
The first one decreases $z$-order of $F_{z^s w^r}$ for the first map by $s'$ and increases
it by the same number for another map. Here both degree and weight do not change.
Finally, we can trade the $w$-order of $F_{z^s w^r}$ the same way
along with moving the appropriate number of $Q$-factors to the other parenthesis.
For instance, we can pass from  $(F_{zw^2}\cdot Q_{z^2}\cdot Q_{z^5}, F_{z^3})$ to
$(F_{zw}\cdot Q_{z^5}, F_{z^3w} \cdot Q_{z^2})$,
where the $Q$-factor $Q_{z^2}$ goes to the second map together with 
the extra derivative in $w$,
whereas the $w$-derivative of the first map decreases.
Again, also here both degree and weight stay clearly the same.

Summarizing, we see that, starting from the above extreme terms
and using the moves as described, we can obtain any other term.
Hence every term has weight $\le 2k_0-1$ as desired.
Furthermore, it follows from the proof that any term in \eqref{transfQQ'}
with all sets of maps $A_{\nu_1,\ldots,\nu_a}$ and $B_{\nu_1^j,\ldots,\nu_{a^j}^j}$ 
being empty, appears only once.
\epf

\subsection{Relations between jets of defining functions and of the Segre varieties}\Label{implicit}
We return to the situation, where $M\subset \C^{n+1}$ 
is a real-analytic hypersurface with a reference point that we continue to assume to be $0$.
As before let $\rho(Z,\bar Z)$ be any defining function of $M$
that we regard as a convergent power series in $(Z,\bar Z)$.
We make a choice of holomorphic coordinates $Z=(z,w)\in \C^n\times \C$ 
such that $\rho_w(0)\ne 0$.
We can then apply the implicit function theorem to the complexified equation
$\rho(z,w,\bar\zeta)=0$ for $(z,w,\zeta)\in \C^n\times \C\times \C^{n+1}$ 
and solve it locally for $w$ in the form $w=Q(z,\bar\zeta)$,
where $Q$ is holomorphic in $(z,\bar\zeta)\in \C^n\times \C^{n+1}$ near $0$.
The function $Q$ can be used to parametrize the Segre varieties: 
$S_{\zeta} = \{ (z, Q(z,\bar \zeta)) : z\in \C^n\}$.

Our goal here will be to establish an explicit relation
between the partial derivatives of $\rho$
at $(0,\bar\zeta)\in \C^{n+1}\times \C^{n+1}$ and of $Q$ at 
$(0,\bar\zeta)\in \C^n\times \C^{n+1}$ 
for $\zeta$ varying in $S_0$,
the Segre variety of $0$ associated to $M$.
We keep the notation $Q_{z^k}$ for the $k$th full derivative of $Q$ in $z$
and use the notation $\rho_{z^s w^l}$ analogous to $G_{z^k w^l}$ in \eqref{g-notation}.
That is, each $\rho_{z^s w^l}(0,\bar\zeta)$ is regarded as an 
$s$-linear function $\C^n\times \cdots\times \C^n\to \C$
depending on the parameter $\zeta\in S_0$.

It turns out that the desired relation has a natural tree structure,
for which we now introduce the needed terminology.
Recall that a {\em (directed or rooted) tree}
is a connected directed graph such 
that each vertex has precisely one incoming arrow except the {\em root}
(one designated vertex) that has none.
We consider here a tree $T$ together with a {\em marking} $s$ by nonnegative integers,
i.e.\ a function $s\colon V(T)\to \N$ (with the convention $\N=\{0,1,\ldots\}$),
where $V(T)$ denotes the set of all vertices of the tree $T$.
The marking will correspond to the differentiation order in $z$.
We do not distinguish between isomorphic marked trees,
i.e.\ trees for which there exist bijections between
their vertices respecting the arrows and the markings.
Together with a marking, we use the integer function $l(a)\in\N$, $a\in V(T)$, 
with $l(a)$ being the number of all outgoing arrows from $a$.
Clearly $l(a)$ depends only on the tree structure (and not on the marking).

\bp\Label{def-segre}
The derivatives of $Q$ and $\rho$ are related by the formula
\begin{equation}\Label{rho-Q}
Q_{z^k}(0,\bar \zeta) = \sum_{T,s} \prod_{a\in V(T)}
\frac{\rho_{z^{s(a)} w^{l(a)}}(0,\bar \zeta)} {-\rho_w(0,\bar \zeta)}, \quad \zeta\in S_0,
\end{equation}
where the product of the multilinear functions is understood in the sense of \eqref{multiplication}
and the summation is taken over the set of all possible finite trees $T$ and their markings $s$ satisfying
\begin{equation}\Label{mark-cond}
2s(a) + l(a) \ge 2\;\; \forall a\in V(T), \quad\sum _{a\in V(T)} s(a) = k.
\end{equation}
\ep

Note that the first condition in \eqref{mark-cond}
eliminates precisely the pairs $(s(a),l(a))$ equal to $(0,0)$ or $(0,1)$.
In particular, the derivative $\rho_w$ in \eqref{rho-Q}
appears only in the denominator.
Note also that both conditions \eqref{mark-cond} together 
force the sum in \eqref{rho-Q} to be finite.
Indeed, summing the inequality in \eqref{mark-cond}
for all vertices $a\in V(T)$ and using the second condition yields 
\begin{equation}\Label{sum}
2k + \sum_a l(a) \ge 2|T|, 
\end{equation}
where $|T|$ stands for the total number of vertices.
Since every vertex has precisely one incoming arrow except the root,
we have $\sum_a l(a)=|T|-1$ by definition of $l(a)$.
Substituting into \eqref{sum} we obtain an estimate on the number of vertices:
\begin{equation}\Label{est}
|T|\le 2k-1.
\end{equation}
Since the number of trees with given number of vertices is finite
and also the number of markings is finite in view of the second condition in \eqref{mark-cond},
we conclude that the sum in \eqref{rho-Q} is finite as claimed.

\br\Label{main-der}
In the sum on the right-hand side of \eqref{rho-Q},
there is precisely one term containing the derivative $\rho_{z^k}$,
namely $\frac{\rho_{z^k}(0,\bar\zeta)}{-\rho_w(0,\bar\zeta)}$,
corresponding to the tree with single vertex $a_0$ and the marking $s(a_0)=k$.
Any other derivative $\rho_{z^s w^l}$ that appears in \eqref{rho-Q},
satisfies $s+l\le k$ and $s<k$.
Indeed, any derivative $\rho_{z^s w^l}$ appears at a vertex $a_0\in T$
with $l$ outgoing arrows. Each outgoing arrow leads, after following a number or arrows, 
to at least one vertex with no further outgoing arrows
(hence corresponding to a derivative $\rho_{z^t}$ with $t\ge 1$). 
Thus we have the vertex $a_0$ with $s(a_0)=s$ and $l$
other vertices $a_1,\ldots,a_l$ with $s(a_j)\ge 1$ for all $j=1,\ldots,l$.
Therefore $\sum_{a\in V(T)} s(a) \ge s+l$ and hence $s+l\le k$
in view of \eqref{mark-cond}. For $s=k$, it must follow that $l=0$
and $s(a)=0$ for any $a\ne a_0\in V(T)$.
The inequality in \eqref{mark-cond} implies $l(a)\ge 2$ for any $a\ne a_0$,
hence any other vertex has at least two outgoing arrows.
But we have seen that each arrow leads to a vertex $a$ with $s(a)\ge 1$.
Hence this is only possible for the tree with the single vertex $a_0$,
proving the claim.
\er

\bpf[Proof of Proposition~\ref{def-segre}]
We shall obtain the formula \eqref{rho-Q}
by differentiating the identity
\begin{equation}\Label{main-id}
\rho(z,Q(z,\bar\zeta),\bar\zeta) = 0
\end{equation}
at $z=0$ and using the induction on $k$. 
Recall that $Q(0,\bar\zeta)=0$ for $\zeta\in S_0$.

For $k=1$, we have 
\begin{equation}\Label{1st-order}
\rho_z(0,\bar\zeta) + \rho_w(0,\bar\zeta)\, Q_z(0,\bar\zeta) =0,
\end{equation}
implying the desired formula in this case,
where the only possible tree $T$ has one vertex $a_0$ and the only possible marking is $s(a_0)=1$.

We now assume the formula for all $k<k_0$ and differentiate \eqref{main-id} $k_0$ times in $z$
at $z=0$ and $\zeta \in S_0$. All derivatives will be understood evaluated at $(0,\bar\zeta)$
as in \eqref{rho-Q} for the rest of the proof and for brevity we shall omit the argument
 $(0,\bar\zeta)$.
With this convention in mind, we obtain:
\begin{equation}\Label{sum-work}
\sum \rho_{z^r w^h} \cdot Q_{z^{k_1}} \cdot \ldots \cdot Q_{z^{k_h}} = 0,
\end{equation}
where the summation is taken over all indices $r,h\in \N$, and for each $h$, all unordered sets of $h$ indices
$k_1,\ldots,k_h\in \N$, satisfying $r+k_1+\cdots + k_h = k_0$.
Note that we continue using the dot for the multiplication defined in \eqref{multiplication}.
The sum \eqref{sum-work} contains precisely one term with $Q_{z^{k_0}}$, namely $\rho_w\cdot Q_{z^{k_0}}$
(which is also the only term with $\rho_w$),
whereas all other derivatives of $Q$ have lower order.
Hence we can solve \eqref{sum-work} for $Q_{z^{k_0}}$ in the form
\begin{equation}\Label{}
Q_{z^{k_0}} = \sum \frac{\rho_{z^r w^h}}{-\rho_w} \cdot Q_{z^{k_1}} \cdot \ldots \cdot Q_{z^{k_h}},
\end{equation}
where now we have the additional restriction $k_j < k_0$ in the sum 
and no factor $\rho_w$ appears in the numerator on the right.
Hence we can use our induction hypothesis and replace each derivative $Q_{z^{k_j}}$
by the right-hand side of \eqref{rho-Q} corresponding to $k=k_j$:
\begin{equation}\Label{product}
Q_{z^{k_0}} = \sum \frac{\rho_{z^r w^h}}{-\rho_w} 
\prod_{(a_1,\ldots,a_h)\in V(T_1)\times \cdots \times V(T_h)} 
\frac{\rho_{z^{s_1(a_1)} w^{l_1(a_1)}}} {-\rho_w}
\cdot \ldots \cdot \frac{\rho_{z^{s_h(a_h)} w^{l_h(a_h)}}} {-\rho_w},
\end{equation}
where the summation is taken over all choices of $h$ trees $T_1,\ldots, T_h$
with markings $s_1,\ldots,s_h$, satisfying 
\begin{equation}\Label{}
2s_j(a_j)+l_j(a_j)\ge 2\;\; \forall a_j\in V(T_j), \quad \sum_{j\in V(T_j)} s_j(a_j) = k_j.
\end{equation}

We now claim that each term in the sum \eqref{product} appears precisely once on the right-hand
side of \eqref{rho-Q} with $k$ replaced by $k_0$.
To show this, we construct for each term a new tree $T$ with marking $s$ as follows.
The vertex set $V(T)$ is the disjoint union of $V(T_1),\ldots, V(T_h)$, and one more vertex $a_0$
that will become the root of $T$. We keep all the arrows within each $T_j$
and add $h$ arrows from $a_0$ to the root of each tree $T_j$.
Finally we keep the marking for each tree $T_j$ and define $s(a_0):=r$ for the root.
It is easy to see that $T$ is again a directed tree and $s$ is a marking satisfying
\eqref{mark-cond} with $k=k_0$. The pair $(T,s)$ constructed this way, yields
precisely the same term in the sum \eqref{rho-Q} as the one we started with.
Vice versa, given a term in \eqref{rho-Q} with $T$ and $s$, we 
can remove the root $a_0\in V(T)$ with its outgoing arrows 
and obtain a finite collection of marked trees $T_1,\ldots,T_h$.
Setting $r:=s(a_0)$, we obtain precisely the same term in \eqref{product}.
Thus we have a one-to-one correspondence between the terms
and hence \eqref{product} implies the desired formula \eqref{rho-Q},
proving it for $k=k_0$.
\epf

\section{Applications to embeddings of hypersurfaces}\Label{appl-embeddings}
\subsection{Linear dependence of partial derivatives}
We now return to our discussion of holomorphic embeddings.
Let $M'\subset\C^{n+m+1}$ be a real 
hyperquadric with a reference point that we shall assume to be the origin $0\in M'$
and denote by $S'_\zeta$ the associated Segre variety of $\zeta\in\C^{n+m+1}$ 
(see \S\ref{implicit}).
Then it follows directly from the definition
that all varieties $S'_\zeta$ are {\em hyperplanes}.
This simple observation will be important in the sequel.

We next  consider a real-analytic submanifold $M\subset M'$ through $0$,
which is generic in a suitable complex submanifold $V\subset \C^{n+m+1}$,
i.e.\ $M\subset V$ and $T_x M + iT_x M = T_x V$ whenever $x\in M$.
The manifold $V$ is also called the {\em intrinsic complexification} of $M$.
Denote by $d$ the real codimension of $M$ in $V$
(which coincides with the CR-codimension of $M$) and set $n:=\dim_\C V - d$ 
(which coincides with the CR-dimension of $M$).
Then the Segre varieties $S_\zeta$ associated to $M$ are $n$-dimensional complex submanifolds of $V$
defined for $\zeta\in V$ near $0$.
We choose complex-linear coordinates $(z,w)\in \C^n\times \C^{m+1}$
vanishing at $0$ such that $S_{\zeta}$ is given by $w=Q(z,\bar\zeta)$,
where $Q$ is a holomorphic function defined in a neighborhood of $0$ in $\C^n\times \1V$, where $\1V$ denotes the conjugate submanifold.
There will be a priori no relation between these coordinates
and those, where $M'$ has the form \eqref{quadric}.
However, we shall only consider linear changes of coordinates for $M'$,
and hence
the property for the Segre varieties $S'_\zeta$ to be hyperplanes remains unchanged.
In the sequel, by the rank of a set of vectors we shall mean
the dimension of their span.

\bl\Label{q'}
Let $M'\subset\C^{n+m+1}$ be a real hyperquadric through $0$
(given by \eqref{quadric} in some linear coordinates) such that 
\begin{equation}\Label{trans-m'}
e:=(0,\ldots,0,1)\notin T_0^c M'.
\end{equation}
Let $M\subset M'\cap V$ be a real-analytic submanifold through $0$ as above, 
whose Segre varieties $S_{\zeta}$ are given by $w=Q(z,\bar\zeta)$, for $(z,w)\in \C^n\times \C^{m+1}$
and $\zeta$ in the intrinsic complexification $V$ of $M$.
Then for any $m+1$ multiindices $\a_1,\ldots,\a_{m+1}\in \N^n$ with $|\a_j|\ge 2$,
the corresponding partial derivatives $Q_{z^{\a_j}}(0,\bar\zeta)$, $j=1,\ldots,m+1$,
are linearly dependent in $\C^{m+1}$ for each $\zeta\in S_0$.
Furthermore, for any given irreducible complex-analytic subvariety $\6S\subset S_0$ passing through $0$,
set
\begin{equation}\Label{max}
r:= \max_{\zeta\in \6S}\rk \{Q_{z^{\a_j}}(0,\bar\zeta) : 1\le j\le m+1 \}\le m.
\end{equation}
Then the first $m$ coordinates of $\C^{m+1}$ can be reordered such that, if 
$\pi \colon \C^{m+1}\to \C^r\times\{0\}\subset \C^r\times\C^{m+1-r}$ stands for the projection
to the first $r$ coordinates, then
\begin{equation}\Label{proj-rank1}
\max_{\zeta\in \6S}\rk \{\pi(Q_{z^{\a_j}}(0,\bar\zeta)):1\le j\le m+1 \} = r.
\end{equation}
\el

\bpf
Recall that $M\subset M'$ implies $S_\zeta\subset S'_\zeta$ for the corresponding Segre varieties (see e.g.\ \cite{W77}).
Consider the parametrization maps $z\mapsto v(z,\bar\zeta):=(z,Q(z,\bar\zeta))$ of the Segre varieties 
$S_{\zeta}$ associated to $M$.
Since for $\zeta\in S_0\subset S'_0$, we have $S_{\zeta}\subset S'_{\zeta}$  and the latter variety is a hyperplane, the derivatives
$$v_{z_1}(0,\bar\zeta),\ldots,v_{z_n}(0,\bar\zeta),
v_{z^{\a_1}}(0,\bar\zeta),\ldots,v_{z^{\a_{m+1}}}(0,\bar\zeta) \in \C^{n+m+1}$$
are also contained in a hyperplane $\Pi$ in $\C^{n+m+1}$ for each $\zeta\in S_0$
(with $\Pi$ depending on $\zeta$).
Since $|\a_j|\ge 2$, we have 
\begin{equation}\Label{vectors}
v_{z^{\a_j}}(0,\bar\zeta)=(0,Q_{z^{\a_j}}(0,\bar\zeta))\in \C^n\times \C^{m+1},
\quad 1\le j\le m+1,
\end{equation}
and therefore these vectors are contained in $\Pi \cap (\{0\}\times \C^{m+1})$.
Since $\Pi= S'_0=T_0^cM'$ for $\zeta=0$, we have $e\notin\Pi$ for $\zeta$ near $0$ in view of \eqref{trans-m'}.
Restricting to a possibly smaller neighborhood of $0$, we may assume that $e\notin\Pi$ holds for all $\zeta$.
Here we use the irreducibility assumption on $\6S$ and its consequence that 
the ranks in \eqref{max} and \eqref{proj-rank1} do not change after restricting
$\zeta$ to any smaller neighborhood of $0$.
Hence $\Pi \cap (\{0\}\times \C^{m+1})$ is a proper hyperplane in $\{0\}\times \C^{m+1}$.
Thus the vectors \eqref{vectors} are linearly dependent. 
Furthermore, since $e$ is not contained in the span of the vectors \eqref{vectors},
the dimension of this span remains unchanged after projecting to the first $m$ coordinates of 
the space $\C^{m+1}$.
Then we can reorder the coordinates of $\C^m\times \{0\}\subset \C^{m+1}$ 
and consider the standard projection $\pi\colon \C^{m+1}\to \C^r\times\{0\}$ such that
\eqref{proj-rank1} holds.
\epf

\subsection{Polynomial relations for the partial derivatives of $Q$}
We now return to the original situation, where $M$ 
is a real-analytic hypersurface in $\C^{n+1}$ with reference point $p\in M$.
As before we choose local holomorphic coordinates $Z=(z,w)\in \C^n\times \C$
vanishing at $p$ and a defining function $\rho(Z,\bar Z)$ for $M$ such that $\rho_w(0,0)\ne0$.
As in \S\ref{implicit} we apply the implicit function theorem to the complexified equation
$\rho(z,w,\bar\zeta)=0$ and solve it locally for $w$ in the form $w=Q(z,\bar\zeta)$,
where $Q$ is a holomorphic function in $(z,\bar\zeta)\in \C^n\times \C^{n+1}$ near $0$
that can be used to parametrize the Segre varieties: 
$S_{\zeta} = \{ (z, Q(z,\bar \zeta)) : z\in \C^n\}$.
We continue to use the weights of the partial derivatives of $Q$ given by 
\begin{equation}\Label{wt3}
\wt Q_{z^{\a}} := 2|\a|-1
\end{equation}
as in \eqref{wt1}.

\bp\Label{Q-relation}
Let $M\subset \C^{n+1}$ be a real-analytic hypersurface through $0$,
which is transversally holomorphically embeddable into a real hyperquadric in $\C^{n+m+1}$.
Then for every irreducible complex-analytic subvariety $\6S\subset S_0$ passing through $0$
and every set of $m+1$ multiindices $\a_j$, $|\a_j|\ge 2$, $j=1,\ldots,m+1$,
there exists an integer $k$ with $K:=\{j : |\a_j| = k\}\ne\emptyset$
such that the partial derivatives of $Q$ satisfy a relation of the form
\begin{equation}\Label{relate1}
\sum_{j\in K} P_j(Q_{z^\b}(0,\bar\zeta))\, Q_{z^{\a_j}}(0,\bar\zeta) = 
R(Q_{z^\b}(0,\bar\zeta)), \quad \zeta\in \6S,
\end{equation}
where $P_j(Q_{z^\b}(0,\bar\zeta))$ and $R(Q_{z^\b}(0,\bar\zeta))$ are some polynomials
in the partial derivatives $Q_{z^\b}$ of lower order (i.e.\ $|\b|<k$),
having weights 
\begin{equation}\Label{est1}
\wt P_j \le (2k-2)(l-1), \quad \wt R \le (2k-2)l+1,
\end{equation}
and not all $P_j(Q_{z^\b}(0,\bar\zeta))$ identically vanish in $\zeta\in \6S$,
where $l\ge1$ is the number of all $j$'s with $|\a_j|\le k$.
\ep

\bpf
Denote by $M'$ a hyperquadric in $\C^{n+m+1}$, where $M$ can be embedded,
and let $H$ be any embedding.
Without loss generality, $M'$ passes through $0$ such that \eqref{trans-m'} holds,
where we use the transversality of the embedding.
We write $H=(F,G)$ as in \eqref{split}.
By a linear change of coordinates in $\C^{n+m+1}$,
we can achieve in addition the normalization assumptions \eqref{normalization} as well as
\begin{equation}\Label{addnorm}
G_w(0)=(0,\ldots,0,1).
\end{equation}
Then we are in the setting of Proposition~\ref{transformation},
where the complex hypersurface $S$ is any Segre variety $S_{\zeta}$ of $M$,
given by $w=Q(z,\bar\zeta)$, $\zeta\in S_0$,
and $S'=H(S)$ is its image in $\C^{n+m+1}$ given by $w'=Q'(z',\bar\zeta)$,
where $Q'$ is an appropriate holomorphic $\C^{m+1}$-valued function in $z'$
with parameter $\zeta$.
The relation between the full higher order derivatives of 
$Q$ and $Q'$ (with respect to $z$ and $z'$ respectively) 
at $0$ is given by \eqref{transfQQ'}.
Since each term's weight on the right-hand side of \eqref{transfQQ'}
does not exceed $2k-1$, only one term can appear with $Q_{z^k}(0,\zeta)$,
namely $G_w(0) Q_{z^k}(0,\zeta)$.
Hence we have 
\begin{equation}\Label{za}
Q'_{z'^\a}(0,\bar\zeta) = G_w(0) Q_{z^\a}(0,\bar\zeta) + P^\a(Q_{z^\b}(0,\bar\zeta)),
\end{equation}
where $P^\a(Q_{z^\b}(0,\bar\zeta))$ is a polynomial in the lower order derivatives $Q_{z^\b}(0,\bar\zeta)$,
$|\b|<|\a|$, with $\wt P^\a \le \wt Q_{z^\a} = 2|\a|-1$.
Moreover, since the derivative $G_w(0)$ satisfies \eqref{addnorm}
and any other derivative $G_{z^s w^l}(0)$ is of positive weight,
we can rewrite \eqref{za} with improved weight estimates as
\begin{equation}\Label{za1}
Q'_{z'^\a}(0,\bar\zeta) = \big(0, Q_{z^\a}(0,\bar\zeta)\big) + 
\big(R^\a(Q_{z^\b}(0,\bar\zeta)), T^\a(Q_{z^\b}(0,\bar\zeta))\big)
\in \C^m\times \C
\end{equation}
with $R^\a$ and $T^\a$ being polynomials of weights 
\begin{equation}\Label{rt-weight}
\wt R^\a \le 2|\a|-2, \quad \wt T^\a \le 2|\a|-1.
\end{equation}

We next apply Lemma~\ref{q'} to the partial derivatives of $Q'$ in $z'$
corresponding to the given multiindices $\a_1,\ldots,\a_{m+1}$.
It follows that there exists an integer $r\le m$
and one can reorder the first $m$ coordinates of $\C^{m+1}$ such that, if 
$\pi \colon \C^{m+1}\to \C^r\times\{0\}$ is the projection to the first $r$ coordinates, then
\begin{equation}\Label{proj-rank}
\max_{\zeta\in \6S}\rk \{Q'_{z'^{\a_j}}(0,\bar\zeta):1\le j\le m+1 \} =
\max_{\zeta\in \6S}\rk \{\pi(Q'_{z'^{\a_j}}(0,\bar\zeta)):1\le j\le m+1 \} = r.
\end{equation}
Without loss of generality, we may assume that the multiindices $\a_j$ are ordered
such that $|\a_1|\le \ldots \le|\a_{m+1}|$.
We claim that an integer $1\le j_0\le m+1$ can be chosen
such that 
\begin{equation}\Label{proj-rank1'}
\max_{\zeta\in \6S}\rk \{\pi(Q'_{z'^{\a_j}}(0,\bar\zeta)):1\le j< j_0 \} =
\max_{\zeta\in \6S}\rk \{\pi(Q'_{z'^{\a_j}}(0,\bar\zeta)):1\le j\le j_0 \} = j_0-1.
\end{equation}
Indeed, denote by $r(j_0)$ the left-hand side of \eqref{proj-rank1'}.
Then $r(j_0)$ is an increasing integer function of $j_0$ with $r(1)=0$ and $r(m+2)\le r\le m$.
Then there must exist $j_0$ with $r(j_0)=r(j_0+1)$ and it suffices
to take the minimum $j_0$ with this property to prove the claim.

We now consider the $j_0\times j_0$ matrix $(Q'^h_{z'^{\a_{j}}}(0,\bar\zeta))$
with $1\le j\le j_0$ and either $1\le h\le j_0-1$ or $h=m+1$. 
Then \eqref{proj-rank1'} implies that, after a suitable permutation of the coordinates in $\C^r$,
the determinant of this matrix identically vanishes, 
whereas the leading $(j_0-1)\times (j_0-1)$ minor corresponding to $1\le j,h\le j_0-1$,
does not identically vanish. In view of \eqref{za1},
the first condition yields an identity of the form \eqref{relate1}
and the second --- the nonvanishing of the coefficient $P_{j_0}(Q_{z^\b}(0,\bar\zeta))$
in front of $Q_{z^{\a_{j_0}}}(0,\zeta)$.
Finally, the desired weight estimates follow from \eqref{rt-weight}.
\epf

\section{Polynomial relations for the derivatives of the defining functions}\Label{poly-relations}
We prove here a stronger version of Theorem~\ref{main},
where we replace the Segre variety $S_p$ of the reference point $p$ 
with any irreducible complex-analytic subvariety of $S_p$.
As before in \S\ref{main-sect} we write $\rho(Z,\bar Z)$ for a real-analytic defining
function of a hypersurface $M\subset\C^{n+1}$ satisfying $\rho_w(p,\bar p)\ne 0$
for some fixed holomorphic coordinates $Z=(z,w)\in \C^n\times \C$
and keep the weights 
$\wt \rho_{z^\a w^s}:= 2|\a|+s-1$
as in \eqref{wt40}.
Recall that the derivative $\rho_w$ is the only one with weight $0$
and we count this derivative separately, denoting by $\deg_{\rho_w} P$
the degree of the polynomial $P$ in the variable $\rho_w$.
We now have the following stronger version of Theorem~\ref{main}:

\bt\Label{main1}
Let $M\subset \C^{n+1}$ be a real-analytic hypersurface through $p$,
which is transversally holomorphically embeddable into a hyperquadric in $\C^{n+m+1}$.
Then for every irreducible complex-analytic subvariety $\6S\subset S_p$ passing through $p$
and every set of $m+1$ multiindices $\a_j$, $|\a_j|\ge 2$, $j=1,\ldots,m+1$,
there exists an integer $k$ with $K:=\{j : |\a_j| = k\}\ne\emptyset$,
such that the partial derivatives of $\rho$ satisfy a relation of the form
\begin{equation}\Label{relate3}
\sum_{j\in K} P_j(\rho_{z^\b w^s}(p,\bar\zeta))\, \rho_{z^{\a_j}}(p,\bar\zeta) = 
R(\rho_{z^\b w^s}(p,\bar\zeta)), 
\quad\zeta\in \6S,
\end{equation}
where $P_j(\rho_{z^\b w^s}(p,\bar\zeta))$ and $R(\rho_{z^\b w^s}(p,\bar\zeta))$ are some polynomials
in the partial derivatives $\rho_{z^\b w^s}(p,\bar\zeta)$ with $|\b|+s\le k$, $|\b|<k$,
and not all $P_j(\rho_{z^\b w^s}(p,\bar\zeta))$ identically vanish in $\zeta\in \6S$.
Moreover, $P_j$ and $R$ can be chosen satisfying in addition the following weight and degree estimates:
\begin{equation}\Label{wt-est}
\begin{array}{lll}
&\wt P_j\le (2k-2)(l-1),\quad &\wt R\le (2k-2)l+1, \\
&\deg_{\rho_w} P_j \le (2k-2)l, \quad &\deg_{\rho_w} R \le (2k-2)l+1,
\end{array}
\end{equation}
where $l$ is the number of all $j$ with $|\a_j|\le k$.
\et

Theorem~\ref{main} corresponds to the special case of Theorem~\ref{main1} with $\6S=S_p$.
Note that the general case of $\6S\subset S_p$ does not follow from that of $\6S=S_p$
by restriction, because all $P_j(\rho_{z^\b w^s}(p,\bar\zeta))$ obtained from Theorem~\ref{main} may identically vanish on the given subvariety $\6S$ even if they don't on $S_p$.
The refined version in the form of Theorem~\ref{main1}
(in fact its proof) will be used in \S\ref{lower} to obtain lower order obstructions
to embeddability than those provided by Theorem~\ref{main}.

\bpf
The proof follows from Propositions~\ref{Q-relation} and \ref{def-segre}.
Indeed, by Proposition~\ref{Q-relation}, we have the relation \eqref{relate1}.
Furthermore, by Proposition~\ref{def-segre}, we can express each derivative of $Q$
by the appropriate expression in the derivatives of $\rho$ according to the formula \eqref{rho-Q}
and substitute them into \eqref{relate1}. 
Multiplying by a suitable power of $\rho_w$, we obtain a polynomial relation \eqref{relate3}.
In view of Remark~\ref{main-der},
the expression substituting for each derivative $Q_{z^\a}$,
contains the derivative $\rho_{z^\a}$ with factor $1\over -\rho_w$ and besides
only the derivatives $\rho_{z^\b w^s}$ with $|\b|+s\le |\a|\le k$ and $|\b|<|\a|\le k$.
Hence the nonvanishing property for the polynomial coefficient in \eqref{relate1} in front of 
some $Q_{z^{\a_j}}$ implies the nonvanishing of the corresponding coefficient
in \eqref{relate3} in front of $\rho_{z^{\a_j}}$.

It remains to show the estimates \eqref{wt-est}.
According to our construction, each derivative $Q_{z^\a}$
is replaced by a sum of terms, each being a product of the derivatives $\rho_{z^\b w^l}$
corresponding to a marked tree $T$ in the formula \eqref{rho-Q},
i.e.\ $|\b|=s(a)$ and $l=l(a)$ for $a\in V(T)$.
Summing the weights \eqref{wt40} for all vertices of $T$,
we obtain the total weight equal to 
\begin{equation}\Label{tot-wt}
2\sum s(a) + \sum l(a) - |T|,
\end{equation}
where $|T|$ stands for the total number of vertices as before.
We have $\sum s(a) = |\a|$ in view of \eqref{mark-cond}.
Recall that $l(a)$ is the number of outgoing arrows from the vertex $a$.
Each vertex has precisely one incoming arrow except the root.
Hence $\sum l(a) = |T| - 1$.
Substituting into \eqref{tot-wt} we obtain that the total
weight of a term replacing $Q_{z^\a}$ is $2|\a|-1$,
which is precisely $\wt Q_{z^\a}$.
Hence our substitution will not change the weights, 
proving the estimates in the first line of \eqref{wt-est}.

To estimate the degree in $\rho_w$, observe that a term substituting for each derivative $Q_{z^\a}$ 
in \eqref{relate1} consists of at most $2|\a|-1 = \wt Q_{z^\a}$ factors $\rho_{z^\b w^l}$
in view of the estimate \eqref{est} for the number of all vertices.
Thus the power of $\rho_w$ in the denominator of a term does not exceed the total weight.
The maximal weight of a term in \eqref{relate1} is $(2k-2)l+1$,
hence the power of $\rho_w$ needed to eliminate the denominators is at most $(2k-2)l+1$.
This proves the estimates in the second line of \eqref{wt-est}.
\epf

\section{Invariants attached to real hypersurfaces}\Label{invariants}

Inspired by Proposition~\ref{Q-relation},
we introduce here series of invariants attached to a germ $(M,p)$ of a real-analytic
hypersurface in $\C^{n+1}$ that provide bounds on possible dimension of a hyperquadric,
where $(M,p)$ can be (transversally) embedded.
As before we choose local holomorphic coordinates $(z,w)\in\C^n\times\C$ near $p$, vanishing at $p$,
such that $M$ is given by $w=Q(z,\bar z,\bar w)$ near $p$ 
with $Q$ being a uniquely determined holomorphic function in its arguments 
$(z,\chi,\tau)\in \C^n\times\C^n\times\C$, defined in a neighborhood of $0$. 
We write $\zeta = (\chi,\tau)\in \C^n\times \C$.

In our first sequence of invariants $r_k(M,p)$
we look for possible relations of the form \eqref{relate1},
ignoring the estimates \eqref{est1}.
More precisely, for every integer $k\ge 2$, define $r_k(M,p)$
to be the maximal number $m$ of the partial derivatives 
$Q_{z^{\a_1}}(0,\zeta),\ldots,Q_{z^{\a_m}}(0,\zeta)$ of order $k$ satisfying no relation of the form
\begin{equation}\Label{relate0'}
\sum_{j=1}^{m} P_j(Q_{z^\b}(0,\zeta))\, Q_{z^{\a_j}}(0,\zeta) = 
R(Q_{z^\b}(0,\zeta)), 
\quad\zeta\in S_0,
\end{equation}
where $P_j(Q_{z^\b}(0,\zeta))$ and $R(Q_{z^\b}(0,\zeta))$ are polynomials
in the lower order partial derivatives $Q_{z^\b}(0,\zeta)$, $|\b|<k$,
and 
$$(P_1(Q_{z^\b}(0,\zeta)),\ldots,P_m(Q_{z^\b}(0,\zeta)))\not\equiv 0.$$

It follows from the transformation rule for the derivatives $Q_{z^\a}(0,\zeta)$
(cf.\ Proposition~\ref{transformation})
that the integers $r_k(M,p)$ so defined, depend only 
on $M$ and $p$ but not on the choice of coordinates $(z,w)$
and hence are biholomorphic invariants of $(M,p)$.
Indeed, the derivatives $Q_{z'^\a}(0,\zeta)$ with $|\a|=k$ in a new coordinate system 
$(z',w')$ are expressed as linear combinations of $Q_{z^\b}(0,\zeta)$ and $1$
with coefficients in the field $\6R$ of all
rational functions in the lower order derivatives $Q_{z^\gamma}(0,\zeta)$, $|\gamma|<k$.
On the other hand, $r_k(M,p)$ can be interpreted
as the dimension of the span of all functions $Q_{z^\a}(0,\zeta)$, $|\a|=k$, 
together with the function $1$, over the field $\6R$.

In our second series, we refine the invariants $r_k(M,p)$
by adding the weight estimates \eqref{est1} to consideration.
We fix some coordinates $(z,w)$ as before and 
define the integers $\3r_k(M,p)$, $k\ge 2$, inductively as follows. 
Assuming that $\3r_k(M,p)$ are defined for $k<k_0$, define $\3r_{k_0}(M,p)$
to be the maximal number $m$ of the partial derivatives 
$Q_{z^{\a_1}}(0,\zeta),\ldots,Q_{z^{\a_m}}(0,\zeta)$ of order $k_0$ 
satisfying no relation of the form
\eqref{relate0'} as above with the additional restriction that
\begin{equation}\Label{}
\wt P_j \le (2k_0-2)\big(\sum_{k<k_0} \3r_k(M,p) + m - 1\big), \quad
\wt P_j \le (2k_0-2)\big( \sum_{k<k_0} \3r_k(M,p) + m \big) +1.
\end{equation}
Analysing the transformation rule given by Proposition~\ref{transformation}
in case of mappings between equal dimension spaces,
we conclude that the integers $\3r_k(M,p)$ remain invariant
under coordinate changes given by $H=(F,G)$ satisfying \eqref{normalization}.
On the other hand, $\3r_k(M,p)$ may potentially change
under the linear coordinate transformations, where the corresponding
change of the derivatives $Q_{z^{\a}}(0,\zeta)$ is rational
rather than polynomial. Thus, in order to obtain an invariant,
we define $\2r_k(M,p)$ to be the minimum of $\3r_k(M,p)$
taken over all possible linear changes of coordinates.

As an immediate consequence of Proposition~\ref{Q-relation},
we now obtain the following relations between the
invariants just defined and embeddings into hyperquadrics:

\bc\Label{hyp-embed}
Let $(M,p)$ be a germ of real-analytic hypersurface in $\C^{n+1}$ that is 
transversally holomorphically embeddable into a real hyperquadric in $\C^{n+m+1}$.
Then 
$$\sum_k r_k(M,p)\le\sum_k\2r_k(M,p) \le m.$$ 
In particular, if $\sum_k r_k(M,p)=\infty$ or $\sum_k \2r_k(M,p)=\infty$, then $(M,p)$ is not 
transversally holomorphically embeddable into any real hyperquadric.
\ec

\br
Similarly to $r_k(M,p)$ and $\2r_k(M,p)$
we can also define further invariants
using the identities \eqref{relate0} in Theorem~\ref{main}
instead of \eqref{relate0'}
and the weight estimates \eqref{wt-est0}
instead of \eqref{est1}.
However, it follows from the proof of Theorem~\ref{main}
(in fact from Proposition~\ref{def-segre})
that so defined invariants do not exceed
$r_k(M,p)$ and $\2r_k(M,p)$ respectively
and hence provide a rougher estimate for
the embeddability dimension.
\er

\section{Distinguished submanifolds of the Segre varieties 
and lower order obstructions}\Label{lower}

The theory of Chern and Moser \cite{CM} reveals some special nature 
of real hypersurfaces of low dimension. 
For instance, in case $n\ge 2$, the obstruction for a real hypersurface
in $\C^{n+1}$ to be a hyperquadric is of order $4$, whereas for $n=1$, it is of order $6$.
This phenomenon turns out to arise in a more elaborated form in our case,
where we study obstructions to embeddability into higher dimensional hyperquadrics.

In order to describe this phenomenon we shall introduce
some distinguished families of submanifolds of the Segre varieties.
Throughout this section $M$ will be a real-analytic Levi-nondegenerate hypersurface in $\C^{n+1}$.
Recall that the family of the Segre varieties $S_Z$, $Z\in \C^{n+1}$,
associated to a generic real-analytic CR-submanifold $M\subset\C^{n+1}$, 
is parametrized by the points of the ambient space.
In case $M$ is a Levi-nondegenerate hypersurface (see e.g. \cite{BERbook}
for this and other basic terminology), 
each $S_Z$ is a complex hypersurface
and the map $S_{p}\ni Z\mapsto T_{p} S_Z$ into the corresponding Grassmannian
is of maximal rank by an observation due to Webster \cite{W77}.
We are going to refine this family as follows.
Given any linear subspace $V\subset T_{p} Q_{p} = T_{p}^cM$, define
\begin{equation}\Label{}
S_{p,V}:=\{ Z\in S_{p} : T_{p} S_Z \supset V \}\subset S_{p}.
\end{equation}
It is easy to see that the sets $S_{p,V}$ are local invariants of $M$,
more precisely, a neighborhood of $p$ in $S_{p,V}$
is completely determined by a neighborhood of $p$ in $M$
and is sent to $S_{p,H_*V}$ (as germ at $p$) by any local biholomorphism $H$ of $\C^{n+1}$ preserving the germ $(M,p)$.
Furthermore, since the map $Z\mapsto T_{p} S_Z$ is of the maximal rank $n$ at $p$,
it follows that each $S_{p,V}$ is a complex submanifold of $\C^{n+1}$ through $p$
(in fact, the tangent space $T_p S_{p,V}$ coincides with the orthogonal
complement of $V$ with respect to the Levi form of $M$).

We keep the notation from \S\ref{j-trans}-\ref{implicit}.
In addition to \eqref{normalization} we assume 
\begin{equation}\Label{normG}
G_z(0)=0.
\end{equation}
The reference point $p\in M$ will be assumed to be $0$.
We also consider the standard basis $e_1,\ldots,e_n$ in $\C^n$
given by $e_j=(0,1,0)\in \C^{j-1}\times \C\times \C^{n-j}$.
For a subset $I\subset \{1,\ldots,n\}$, we set
\begin{equation}\Label{}
V^0_I:=\span \{e_j : j\in I\}, \quad V_I:=(V^0_I\times \C) \cap T^c_0M.
\end{equation}
For every such $I$, consider the distinguished submanifold $S_{0,V_I}\subset S_0$.
We also use the notation 
$$\supp \a:= \{j: \a_j\ne 0 \} \subset \{1,\ldots,n\}$$ 
for a multiindex $\a=(\a_1,\ldots,\a_n)\in \N^n$.

\bl\Label{order2}
Let $I\subset \{1,\ldots,n\}$ be a proper subset and
$\a$ be a multiindex with $|\a|=2$ and $\supp \a\subset I$.
Then
\begin{equation}\Label{deg2}
Q'_{z^\a}(0,\bar\zeta) = G_w \left(Q_{z^\a}(0,\bar\zeta) 
+ P_1(Q_{z^\b}(0,\bar\zeta))\right) + P_0,\quad \zeta \in S_{0,V_I},
\end{equation}
where $P_1(Q_{z^\b}(0,\bar\zeta)$ is a polynomial of degree at most $1$
in the derivatives $Q_{z^\b}(0,\bar\zeta)$, $|\b|=1$, and $P_0$ is a constant.
\el

\bpf
We apply the transformation formula \eqref{transfQQ'} for the
$2$nd full derivatives of $Q(z,\bar \zeta)$ with respect to $z$
that we denote by $Q_{z^2}(z,\bar \zeta)$.
We have $Q'_{z^\a}(0,\bar\zeta)=Q_{z^2}(z,\bar \zeta)(v^1,v^2)$ 
for suitable vectors $v^1,v^2$ (from the standard basis $e_1,\ldots,e_n$).
Hence, in order to obtain $Q'_{z^\a}(0,\bar\zeta)$,
we evaluate \eqref{transfQQ'} with $k=2$ at $(v^1,v^2)$.
The key observation is that $\supp \a\subset  I$
implies $Q_{z}(0,\bar \zeta)(v^j)$ being constant in $\bar \zeta$
for $j=1,2$, and $\zeta \in S_{0,V_I}$.
Hence, on the right-hand side of \eqref{transfQQ'},
the only other factors involving $Q_z$
are of the form $(Q_z\circ B_0)$, $(Q_z\circ B_1)(v^j)$, $(Q_z\circ B_2)(v^1,v^2)$
having degrees $0,1,2$ and weights $1,2,3$ respectively,
whereas $Q_{z}(0,\bar \zeta)(v^j)$ has both degree and weight $1$.
According to Proposition~\ref{transformation},
the degree of each term is $2$ and the weight does not exceed $3$.
Hence all terms contain only the first order derivatives $Q_{z^\b}(0,\bar\zeta)$, $|\b|=1$,
except $G_w Q_{z^\a}(0,\bar\zeta)$ (recall that $\wt Q_{z^\a} = 3$), the latter giving the first term 
on the right-hand side of \eqref{deg2}.
Since the total degree must be $2$ and total weight $\le 3$, any other term with $G_w$ has either 
at most one factor $Q_z\circ B_1$ or at most one factor $Q_z\circ B_2$. 
In the first case, the other factor has to be $Q_z$ which is constant, hence the corresponding product is a polynomial 
in $Q_{z^\b}(0,\bar\zeta)$, $|\b|=1$, of degree at most $1$,
contributing to the polynomial $P_1$ in \eqref{deg2}.
In the second case, $\wt Q_z\circ B_2 = 1+2 =3$, hence no other factor can appear
and the result is again a polynomial in $Q_{z^\b}(0,\bar\zeta)$, $|\b|=1$, of degree at most $1$,
contributing to $P_1$.

Finally, we investigate the terms containing the other derivatives $G_{z^kw^l}$
(i.e.\ with $(k,l)\ne (0,1)$).
Some of them contain only constant factors with $Q_z$ 
and hence contribute to $P_0$ in \eqref{deg2}.
Any other term which is not constant, must have a factor $Q_z\circ B_\nu$ and thus is of weight at least $1$.
Therefore we can only have $G_{z^kw^l}$ or $G_{z^kw^l}\circ A_\mu$ with total weight $\le 2$.
Since $(k,l)\ne (0,1)$ and in view of \eqref{normG}, we can only have 
$G_{w^2}$, $G_{w^3}$, $G_{zw}$ or $G_{zw}\circ A_0$.
Since the total degree has to be $2$, each term has
a factor $Q_z\circ B_\nu$ with $\nu\ge 1$, having weight $\ge 2$.
This leaves  a weight at most $1$ for the factor involving $G_{z^kw^l}$, which can only be $G_{w^2}$. But $G_{w^2}$ requires at least two factors involving $Q_z$, each having weight $\ge1$.
The latter makes it impossible to have the total weight not exceeding $3$ and therefore no terms of that kind may occur.
\epf

Using Lemma~\ref{order2} we obtain low degree relations
between second order derivatives:

\bc\Label{Q-relate2}
Let $M\subset \C^{n+1}$ be a real-analytic hypersurface through $0$
given by $w=Q(z,\bar Z)$, $Z=(z,w)\in \C^n\times \C$,
that is transversally holomorphically embeddable into a real hyperquadric in $\C^{n+m+1}$.
Fix a proper subset $I\subset\{1,\ldots,n\}$.
Then for every set of $m+1$ multiindices $\a_j$, $j=1,\ldots,m+1$, 
with $|\a_j|=2$ and $\supp \a_j \subset I$,
there is a relation of the form
\begin{equation}\Label{relate1'}
\sum_{j} \l_j\, Q_{z^{\a_j}}(0,\bar\zeta) = 
R_1(Q_{z^\b}(0,\bar\zeta)), \quad \zeta\in S_{0,V_I},
\end{equation}
where $(\l_1,\ldots,\l_{m+1})\ne 0$ and $R_1(Q_{z^\b}(0,\bar\zeta))$ 
is a polynomial of degree at most $1$ 
in the first order partial derivatives $Q_{z^\b}(0,\bar\zeta)$.
\ec

\bpf
The proof is analogous to that of Proposition~\ref{Q-relation},
where we use Lemma~\ref{order2} instead of Proposition~\ref{transformation}.
The details are left to the reader.
\epf

Using the transformation law between the derivatives of $Q$ and of $\rho$
given by Proposition~\ref{def-segre},
we also obtain a low degree relation between the derivatives of $\rho$:

\bc\Label{rho-relate2}
Let $M\subset \C^{n+1}$ be a real-analytic hypersurface through $0$
given by $\rho(Z,\bar Z)=0$, $\rho_w(0,0)\ne 0$,
that is transversally holomorphically embeddable into a real hyperquadric in $\C^{n+m+1}$.
Fix a proper subset $I\subset\{1,\ldots,n\}$.
Then for any set of $m+1$ multiindices $\a_j$, $j=1,\ldots,m+1$, 
with $|\a_j|=2$ and $\supp \a_j \subset  I$,
there is a relation of the form
\begin{equation}\Label{relate42}
\sum_{j} \l_j\, \frac{\rho_{z^{\a_j}}(0,\bar\zeta)}{\rho_w(0,\bar\zeta)} = 
R_1\left(
\frac{\rho_{z^\b}(0,\bar\zeta)}{\rho_w(0,\bar\zeta)},
\frac{\rho_{z^\b w}(0,\bar\zeta)}{\rho_w(0,\bar\zeta)},
\frac{\rho_{w^2}(0,\bar\zeta)}{\rho_w(0,\bar\zeta)}
\right), \quad \zeta\in S_{0,V_I},
\end{equation}
where $(\l_1,\ldots,\l_{m+1})\ne 0$ and $R_1$ 
is a polynomial of degree at most $1$ in its components
involving $|\b|\le 1$.
\ec

\bpf
As mentioned before, we use the transformation law between the derivatives of $Q$ and $\rho$
given by Proposition~\ref{def-segre}. The relation between 
$Q_{z^\b}(0,\bar\zeta)$ and $\rho_{z^\b}(0,\bar\zeta)$
follows, for instance, from \eqref{1st-order}.
To obtain the formula for $Q_{z^{\a_j}}(0,\bar\zeta)$,
we use \eqref{rho-Q} for $k=2$:
\begin{equation}\Label{k=2}
Q_{z^2}(0,\bar\zeta) = 
\frac{\rho_{z^2}(0,\bar\zeta)}{-\rho_w(0,\bar\zeta)} + 
\frac{\rho_{zw}(0,\bar\zeta)}{-\rho_w(0,\bar\zeta)}
\frac{\rho_{z}(0,\bar\zeta)}{-\rho_w(0,\bar\zeta)}  +
\frac{\rho_{w^2}(0,\bar\zeta)}{-\rho_w(0,\bar\zeta)}
\left(\frac{\rho_{z}(0,\bar\zeta)}{-\rho_w(0,\bar\zeta)}\right)^2.
\end{equation}
The partial derivative $Q_{z^{\a_j}}(0,\bar\zeta)$
is now given by the evaluation of the right-hand side of \eqref{k=2}
at the suitable pair of vectors $(v^1,v^2)$.
By the same key observation as in the proof of Lemma~\ref{order2},
we conclude that $Q_{z}(0,\bar\zeta)(v^j)=0$
and hence the ratio $\frac{\rho_{z}(0,\bar\zeta)}{-\rho_w(0,\bar\zeta)}$ 
is constant in $\zeta\in S_{0,V_I}$.
Using this information and substituting the obtained formulas in \eqref{relate1'}
we come to the desired conclusion.
\epf

Note that in Corollary~\ref{rho-relate2}
there are no restrictions on $\rho$ and the coordinates chosen
other than $\rho_w(0,0)\ne 0$.
However, if we choose $(z,w)$ such that the complex tangent space $T_0^cM$
is given by $w=0$, the conclusion of Corollary~\ref{rho-relate2}
is substantially simplified:

\bc\Label{rho-relate3}
Let $M\subset \C^{n+1}$ be a real-analytic hypersurface through $0$
given by $\rho(Z,\bar Z)=0$ with $\rho_w(0,0)\ne 0$ and $\rho_z(0,0)=0$.
Suppose that $M$ is transversally embeddable into a real hyperquadric in $\C^{n+m+1}$.
Fix a proper subset $I\subset\{1,\ldots,n\}$.
Then for any set of $m+1$ multiindices $\a_j$, $j=1,\ldots,m+1$, 
with $|\a_j|=2$ and $\supp \a_j \subset  I$,
there is a relation of the form
\begin{equation}\Label{relate52}
\sum_{j} \l_j\, \rho_{z^{\a_j}}(0,\bar\zeta) = 
R_1(\rho_{z^\b}(0,\bar\zeta)), \quad \zeta\in S_{0,V_I},
\end{equation}
where $(\l_1,\ldots,\l_{m+1})\ne 0$ and $R_1$ 
is a polynomial of degree at most $1$ in 
the first order derivatives of $\rho$.
\ec

\bpf
Recall from the proof of Corollary~\ref{rho-relate2} that 
the ratio $\frac{\rho_{z}(0,\bar\zeta)}{-\rho_w(0,\bar\zeta)}$ 
is constant in $\zeta\in S_{0,V_I}$.
Since now we assume $\rho_z(0,0)=0$, this ratio is actually zero.
Then \eqref{k=2} is reduced to 
$Q_{z^2}(0,\bar\zeta) = 
\frac{\rho_{z^2}(0,\bar\zeta)}{-\rho_w(0,\bar\zeta)}$.
The rest of the proof is completely analogous to that of Corollary~\ref{rho-relate2}.
\epf

\br
As in Remark~\ref{rem-main}, we consider the special case $m=0$, 
where Corollaries~\ref{Q-relate2} -- \ref{rho-relate3} 
give obstructions preventing $M$ from being (locally) biholomorphically equivalent to 
a hyperquadric
and make a comparison with the Chern-Moser normal form
\begin{equation}\Label{normalform}
M=\big\{\Im w = \sum a_{\a\mu s} z^\a \bar z^\mu (\Re w)^s\big\},
\end{equation}
where, in particular, there are no pure terms $z^\a$ and $(\Re w)^s$
and the Levi form of $M$ at $0$ is given by $\sum \pm |z_j|^2$.
If all the first order derivatives $\rho_{z^\b}(0,\bar\zeta)$ are (affine) linear functions
(which is the case in the Chern-Moser normal form),
Corollary~\ref{rho-relate3} implies that, in case $M$ is equivalent to a hyperquadric, 
all second order derivatives 
$\rho_{z^{\a}}(0,\bar\zeta)$ with $\supp \a \subset I$
are linear in $\zeta\in S_{0,V_I}$.
Since $S_0$ is given by $w=0$, the property $\zeta\in S_{0,V_I}$ means
$\zeta_s=0$ for $s\in I$ (i.e.\ $\zeta$ is orthogonal to $V_I$
with respect to the Levi form).
Now the mentioned linearity of $\rho_{z^{\a}}(0,\bar\zeta)$ for $\zeta\in S_{0,V_I}$
means that $a_{\a\mu 0}=0$ whenever $|\a|=2$, $|\mu|\ge 2$ and 
\begin{equation}\Label{supp-restriction}
\supp \a\cap \supp \mu = \emptyset.
\end{equation}
On the other hand, if $M$ is equivalent to a hyperquadric, 
the Chern-Moser theory implies the vanishing of certain
terms of bidegree $(2,2)$ and $(2,3)$ in $(z,\bar z)$ in the normal form.
One can see that \eqref{supp-restriction} is closely related
to the trace-free parts of the corresponding polynomials (see \cite[p. 233]{CM}).
However, as mentioned before, the actual normal form can be hard to calculate,
whereas Corollaries~\ref{Q-relate2} -- \ref{rho-relate3} can be applied directly
in any given coordinates.
\er

We conclude this section by a series of explicit examples of manifolds $M$
that are not embeddable into hyperquadrics of certain dimensions
by means of low order obstructions.

\be
Consider any submanifold $M\subset\C^{n+1}$ given by
\begin{equation}\Label{form}
\rho:=-\Im w + \sum_{s=1}^n \pm|z_s|^2 + \sum_{|\a|+k,|\b|+l\ge 2} 
\rho_{\a k \b l} z^\a w^k \bar z^\b \bar w^l = 0.
\end{equation}
In fact, any $M$ with nondegenerate Levi form can be written as \eqref{form},
which is a part of the Chern-Moser normalization.
It easily follows that $S_0=\{w=0\}$, $\rho_w(0,\bar\zeta)=\const$ and 
$\rho_z(0,\bar\zeta)$ is linear in $\bar \zeta$.
We choose two sets of different multiindices $\a_1,\ldots,\a_{m+1}$ and 
$\b_1,\ldots,\b_{m+1}$ with $|\a_j|=2$, $|\b_j|\ge 2$, such that
\begin{equation}\Label{ineq}
\big(\bigcup_j \supp \a_j\big) \cap \big(\bigcup_j \supp \b_j\big) = \emptyset
\end{equation}
and consider the determinant $A$ of the matrix $(\rho_{\a_j0\b_k0})_{jk}$.
Then if $A\ne 0$, Corollary~\ref{rho-relate3} implies that 
$M$ is not transversally embeddable into any hyperquadric in $\C^{n+m+1}$.
In particular, if 
$$m+1\le \frac12 \left[\frac{n}2\right] \left(\left[\frac{n}2\right] - 1\right),$$ 
we can always choose $\b_k$ with $|\b_k|=2$ 
and thus have an obstruction of order $4$.
Indeed, given \eqref{ineq}, 
we can split the set $\{1,\ldots,n\}$ into disjoint subsets $I_1$ and $I_2$ with $[n/2]$ elements each
and choose $\a_j,\b_k$ with $|\a_j|=|\b_k|= 2$ such that $\supp \a_j \subset I_1$ and $\supp \b_k \subset I_2$.
\ee

\section{Embeddability of submanifolds of higher codimension}\Label{high}
\subsection{Obstructions to embeddability}

Our goal here will be to extend some of the preceding results
from hypersurfaces to generic submanifolds $M\subset \C^{n+d}$ of arbitrary codimension $d$.

We begin by giving a version of Proposition~\ref{transformation},
where we adopt all the notation from \S\ref{j-trans},
except that we consider a holomorphic embedding 
$H=(F,G)\colon (\C^{n+d},0) \to (\C^n\times \C^{m+d},0)$ (i.e.\ $1$ is replaced with
general codimension $d$),
$S\subset \C^{n+d}$ is a complex submanifold of codimension $d$ through $0$,
and choose the coordinates $(z,w)\in \C^n\times \C^d$.
As before, $S$ and $S'$ are respectively graphs of holomorphic functions
$w=Q(z)$ and $w'=Q'(z')$ near $0$ with $Q(0)=0$, $Q'(0)=0$
and all derivatives of $Q$ and $Q'$ will be assumed taken at $0$.
We continue writing $Q_{z^\a}\in \C^d$ for a partial derivative 
with respect to a multiindex $\a\in \N^n$
and denote by $Q^i_{z^\a}\in \C$ the components for $1\le i\le d$.
Similar notation is used for $Q'$.
We regard the derivative $G_w$ as an $m\times d$ matrix.

\bp\Label{transformation-high}
Under the normalization assumption \eqref{normalization},
the partial derivatives of $Q$ and $Q'$ at $0$ are related by the formula
\begin{equation}\Label{transf-high}
Q'_{z'^\a} = G_w Q_{z^\a} + P_\a(Q^i_{z^\b}),
\end{equation}
where $P_\a(Q^i_{z^\b})$ is a $\C^m$-valued polynomial in the components
of the lower order derivatives $Q^i_{z^\b}$, $|\b|<|\a|$.
\ep

\bpf
The proof follows the line of the proof of Proposition~\ref{transformation},
involving differentiation of \eqref{basic} and using induction on $|\a|$.
It is clear that $G_w Q_{z^\a}$ is the only term on the right-hand side
of \eqref{transf-high} involving derivatives of $Q$ of order $k$.
The remainder is a polynomial in the components of the lower order derivatives.
The details are left to the reader.
\epf

We next give a version of Proposition~\ref{def-segre},
relating the derivatives of $Q$ and $\rho$ in the {\em same} coordinates
$(z,w)\in \C^n\times \C^d$.
This time both $Q(z,\zeta)$, $(z,\zeta)\in \C^d \times \C^{n+d}$, and 
$\rho(z,\zeta)$ are $\C^d$-valued
and the coordinates are chosen such that $\rho_w(0,0)$
is an invertible $d\times d$ matrix.
We write $\rho^i_{z^{\b} w^{\gamma}}\in \C$ for the components
of the partial derivatives corresponding 
to integers $1\le i\le d$ and multiindices $\b\in \N^n$ and $\gamma\in \N^d$.

\bp\Label{def-segre-high}
The derivatives of $Q$ and $\rho$ are related by the formula
\begin{equation}\Label{rho-Q1}
Q_{z^\a}(0,\bar \zeta) = -\rho_w^{-1}(0,\bar \zeta) \rho_{z^\a}(0,\bar \zeta)
+ \frac{R_\a(\rho^i_{z^{\b} w^{\gamma}}(0,\bar \zeta))}
{(\det \rho_w(0,\bar \zeta))^{l_\a}}, 
\quad \zeta\in S_0,
\end{equation}
where $R_\a(\rho^i_{z^\b w^{\gamma}}(0,\bar \zeta))$
is a $\C^d$-valued polynomial in the partial derivative components
$\rho^i_{z^{\b} w^{\gamma}}(0,\bar \zeta)$ with $|\b|+|\gamma|\le |\a|$,
$|\b|<|\a|$, and $l_\a$ is a positive integer.
\ep

\bpf
Here we follow the line of the proof of Proposition~\ref{def-segre}.
As in that proof, we differentiate \eqref{main-id},
this time a vector identity, and subsequently use
induction on $|\a|$ when substituting for the components of
$Q_{z^\b}(0,\bar \zeta)$ with $|\b|<|\a|$.
The details are left to the reader.
\epf

We now turn to a version of Theorem~\ref{main1} for higher codimension.

\bt\Label{high-embed}
Let $M\subset \C^{n+d}$ be a real-analytic generic submanifold through $0$ given by 
$\rho(Z,\bar Z)=0$
with $\rho_w(0,0)$ being invertible $d\times d$ matrix.
Suppose that $M$ is transversally holomorphically embeddable 
into a hyperquadric in $\C^{n+m+d}$.
Then for any set of $m+d$ multiindices $\a_j\in \N^n$, $|\a_j|\ge 2$, $j=1,\ldots,m+d$,
there exist  integers $i_0\in \{1,\ldots,d\}$ and $k$ with $K:=\{j: |\a_j| = k\}\ne \emptyset$
such that the partial derivative components of $\rho$ satisfy a relation of the form
\begin{equation}\Label{relate4}
\sum_{j\in K} P_j(\rho^i_{z^\b w^\gamma}(0,\bar\zeta))\, \rho^{i_0}_{z^{\a_j}}(0,\bar\zeta) = 
R(\rho^i_{z^\b w^\gamma}(0,\bar\zeta)), 
\quad\zeta\in S_0,
\end{equation}
where $P_j(\rho^i_{z^\b w^\gamma}(0,\bar\zeta))$ and $R(\rho^i_{z^\b w^\gamma}(0,\bar\zeta))$ 
are polynomials in the partial derivative components 
$\rho^i_{z^\b w^\gamma}(0,\bar\zeta)$ with $|\b|+|\gamma|\le k$ and either $|\b|<k$ or $i\ne i_0$,
such that not all $P_j(\rho^i_{z^\b w^\gamma}(0,\bar\zeta))$
identically vanish in $\zeta\in S_0$.
\et

\bpf
We follow the strategy of the proof of Theorem~\ref{main1}.
We first establish a version of Proposition~\ref{Q-relation}.
As in the proof of the latter, assuming a transversal embedding
$H=(F,G)$ is given, we can perform a linear change
of the coordinates in the target space
and a possible permutation of the components of $w\in \C^d$ in the source to obtain
$$(0,\ldots,0,1)\notin T_0^c M', \quad
G_w (0) =(0,\id)\colon \C^d \to \C^{m}\times \C^d.$$
Then we apply Proposition~\ref{transformation-high}
(in place of Proposition~\ref{transformation})
to obtain a relation 
\begin{equation}\Label{za2}
Q'_{z'^\a}(0,\bar\zeta) = \big(0, Q_{z^\a}(0,\bar\zeta)\big) + 
\big(R^\a(Q^j_{z^\b}(0,\bar\zeta)), T^\a(Q^j_{z^\b}(0,\bar\zeta))\big)
\in \C^m\times \C^d,
\end{equation}
with $R^\a$ and $T^\a$ being polynomials in the lower order derivatives components
$Q^j_{z^\b}(0,\bar\zeta)$, $|\b|<|\a|$.

As in the proof of Proposition~\ref{Q-relation},
we next apply Lemma~\ref{q'} to the given $m+d$ (instead of $m+1$) multiindices 
$\a_1, \ldots, \a_{m+d}$ (assumed to be ordered as $|\a_1|\le \ldots \le|\a_{m+d}|$)
 to obtain an integer $j_0$ and a $j_0\times j_0$ matrix 
$(Q'^h_{z'^{\a_{j_s}}}(0,\bar\zeta))$ with vanishing determinant,
whose leading $(j_0-1)\times (j_0-1)$ minor does not identically vanish.
Then substituting the right-hand side expressions from \eqref{za2}
for the matrix entries, we obtain  a relation
\begin{equation}\Label{relate2'}
\sum_{j\in K} P_j(Q^i_{z^\b}(0,\bar\zeta))\, Q^d_{z^{\a_j}}(0,\bar\zeta) = 
R(Q^i_{z^\b}(0,\bar\zeta)), \quad \zeta\in \6S,
\end{equation}
where $K:=\{j: |\a_j| = k\}\ne \emptyset$ for $k:= |\a_{j_0}|$ and
$P_j(Q^i_{z^\b}(0,\bar\zeta))$ and $R(Q^i_{z^\b}(0,\bar\zeta))$ are polynomials
in the partial derivatives components $Q^i_{z^\b}(0,\bar\zeta)$ 
with either $|\b|<k$ or $|\b|=k$ and $i\ne d$,
and such that not all $P_j(Q^i_{z^\b}(0,\bar\zeta))$ identically vanish.
Note that we previously made a possible permutation of the components of $w$,
so that the last component $Q^d_{z^{\a_j}}(0,\bar\zeta)$ in \eqref{relate2'}
may actually correspond to another component $Q^{i_0}_{z^{\a_j}}(0,\bar\zeta)$
in the original numeration.

Finally we follow the line of the proof of Theorem~\ref{main1},
where we apply Proposition~\ref{def-segre-high}
instead of Proposition~\ref{def-segre}
to pass from the identity \eqref{relate2'}
to an identity of the form \eqref{relate4} as desired.
\epf

\subsection{Most generic submanifolds of higher codimension are not embeddable}
Our goal here is to use Theorem~\ref{high-embed}
in order to give an affirmative answer to a question by Forstneri\v c \cite{Fo04}.
Informally speaking, this question is whether 
{\em the set of all generic submanifolds of higher codimension,
which are holomorphically embeddable into algebraic strongly pseudoconvex hypersurfaces,
is of the first category}.

To state the question more precisely, let us recall some notation from \cite{Fo04}.
Recall that every germ of a generic real-analytic submanifold $M\subset \C^{n+d}$ of codimension $d$
is biholomorphically equivalent to one of the form
\begin{equation}\Label{germ}
M=\big\{ \Im w =  r (\Re z,\Im z,\Re w)\big\}, \quad (z,w)\in \C^n\times \C^d,
\end{equation}
where 
\begin{equation}\Label{series}
r(x,y,u) = \sum_{\a,\b\in\N^n, \g\in\N^d} c_{\a\b\g} x^\a y^\b u^\g, \quad x,y\in \R^n, u\in \R^d,
\end{equation}
is a $\R^d$-valued convergent power series without constant and linear terms.
Then all convergent power series in \eqref{series}
can be written as $\cup_{t>0} \6R^t$,
where $\6R^t$ is the space of all series \eqref{series} for which the norm
\begin{equation}\Label{}
\|r\|_t := \sum_{\a,\b,\g} |c_{\a\b\g}| t^{|\a|+|\b|+|\g|}
\end{equation}
is finite. Clearly $\6R^t$ is a Banach space with the norm $\|r\|_t$.
Finally recall that a real submanifold $M\subset \C^{n+d}$ is called {\em algebraic}
if it is contained in a real-algebraic variety of $\C^{n+d}\cong \R^{2(n+d)}$ of the same dimension as $M$.

We now state our main result of this section.

\bt\Label{nowhere}
For every $t>0$, the set of all $r\in \6R^t$,
for which the germ $(M,0)$ given by \eqref{germ} is transversally holomorphically embeddable into a hyperquadric
in any dimension, is of the first category (in the Banach space $\6R^t$).
\et

Theorem~\ref{nowhere} answers the above question by Forstneri\v c
in view of the result by Webster \cite{W78} stating
that any Levi-nondegenerate {\em real-algebraic} hypersurface 
is always transversally holomorphically embeddable into a Levi-nondegenerate hyperquadric
(of possibly high dimension depending on the hypersurface).

\bpf[Proof of Theorem~\ref{nowhere}]
We first rewrite the power series $r$ in the complex form:
\begin{equation}\Label{ser1}
r = \sum_{\a,\b\in\N^n, \g\in\N^d} r_{\a\b\g} z^\a \bar z^\b u^\g, \quad z\in \R^n, u\in \R^d.
\end{equation}
Then we can identify the elements of $\6R^t$ with the power series \eqref{ser1}
without constant and linear terms whose coefficients $r_{\a\b\g}$ satisfy the reality condition
$\1{r_{\a\b\g}}= r_{\b\a\g}$.
As the next step we eliminate all pure terms $r_{\a00}z^\a$ by subtracting them from $r$.
The corresponding transformation is biholomorphic and hence 
does not change the biholomorphic equivalence class of $(M,0)$.
Denote by $\6R^t_0\subset \6R^t$ the subspace of all series with $r_{\a00}=0$ for all $\a$.
Then it is sufficient to prove the statement for $\6R^t_0$,
i.e.\ to show that the set of all $(M,0)$ corresponding to elements in $\6R^t_0$,
that are transversally embeddable into a hyperquadric, is of the first category in $\6R^t_0$.

We now consider germs $(M,0)$ given by some $r\in \6R^t_0$
that are transversally embeddable into a hyperquadric in $\C^{n+m+d}$ for some fixed $m$.
For every such $(M,0)$, we can apply Theorem~\ref{high-embed} 
using the defining function $\rho$ of $M$ given by
$$\rho(z,w,\bar z,\bar w):= - \Im w + r(\Re z,\Im z,\Re w)$$
and obtain a relation \eqref{relate4}.
Since $r$ has no pure terms with $z^\a$, we have $\rho(z,0,0,0)\equiv 0$, implying that $S_0 = \{w=0\}$.
Then Theorem~\ref{high-embed} yields, in particular, for some $i_0$ and $k$, a polynomial identity
\begin{equation}\Label{relate3'}
\sum_{|\a|=k} P_\a(\rho^i_{z^\b w^\gamma}(0,\bar\chi,0))\, \rho^{i_0}_{z^{\a}}(0,\bar\chi,0) = 
R(\rho^i_{z^\b w^\gamma}(0,\bar\chi,0)), \quad \zeta = (\chi,\tau)\in \C^n\times \C^d,
\end{equation}
with $P_\a$ and $R$ being polynomials as in Theorem~\ref{high-embed}
and not all $P_\a(\rho^i_{z^\b w^\gamma}(0,\bar\chi,0))$ identically vanishing.

In our next step we consider the following standard lexicographic order on the set of all
multiindices $\a=(\a^1,\ldots,\a^n)\in\N^n$. We write $\a < \b$ if either $|\a|<|\b|$
or $|\a|=|\b|$ and for some $1\le s \le n$, $\a^j=\b^j$ for all $j < s$ but $\a^s < \b^s$.
We also write $\a\le\b$ if either $\a <\b$ or $\a =\b$.
Then the set $\N^n$ becomes linearly ordered with the following additive property:
\begin{equation}\Label{additive}
\a_1 \le \b_1, \quad \a_2 \le \b_2 \quad \Longrightarrow \quad \a_1+\a_2 \le \b_1 + \b_2.
\end{equation}

We now fix $i_0$, $k$ and a multiindex $\a_0$ with $|\a_0|=k$
and consider the set of all $r\in \6R^t_0$,
for which a relation \eqref{relate3'} holds with the coefficient
$P_{\a_0}(\rho^i_{z^\b w^\gamma}(0,\bar\chi,0))\not\equiv 0$.
Using the lexicographic order introduced above
we may consider the minimal multiindex 
$\nu_0=\nu_0\big(P_{\a_0}(\rho^i_{z^\b w^\gamma}(0,\bar\chi,0))\big)$ corresponding
to a nonzero monomial in the expansion of $P_{\a_0}(\rho^i_{z^\b w^\gamma}(0,\bar\chi,0))$.
In addition to the previous data, we also fix this minimal multiindex $\nu_0$
as well as the degrees of the polynomials $P_{\a}$ and $R$.
It is clearly sufficient to prove that the set of all $r\in \6R^t_0$
with $\rho$ satisfying \eqref{relate3'} with the above data fixed,
is of the first category.

Going back to \eqref{relate3'}, we 
expand both sides as power series in $\bar\chi$
and obtain recursive relations for the terms of 
$\rho^{i_0}_{z^{\a_0}}(0,\bar\chi,0)$ as follows.
For every multiindex $\mu$, $|\mu|\ge1$, identify the monomials in the expansion with $\bar\chi^{\nu_0+\mu}$.
Then, since $P_{\a_0}(\rho^i_{z^\b w^\gamma}(0,\bar\chi,0))$
contains a nontrivial monomial with $\bar\chi^{\nu_0}$,
we have a nontrivial contribution of the corresponding monomial with $\bar\chi^{\mu}$
in the expansion of $\rho^{i_0}_{z^{\a_0}}(0,\bar\chi,0)$.
(The latter monomial may be assumed nonvanishing, since it vanishes only for a set of $r$'s of the first category.)
Furthermore, it follows from the property of $\nu_0$ and \eqref{additive}
that the contributing multiindices corresponding to all 
other nontrivial monomials in the expansion of $\rho^{i_0}_{z^{\a_0}}(0,\bar\chi,0)$
are smaller than $\mu$.
Thus we can express the coefficient in front of $\bar\chi^{\mu}$ in the expansion of 
$\rho^{i_0}_{z^{\a_0}}(0,\bar\chi,0)$ as a rational function of its other coefficients
corresponding to smaller monomials (with respect to our lexicographic order), the coefficients of other derivatives 
$\rho^{i}_{z^\b w^\g}(0,\bar\chi,0)$ and the coefficients of $P_{\a}$ and $R$.
By induction, we can then express the coefficients of $\bar\chi^{\mu}$
as a rational function only of the coefficients of other derivatives and the polynomials 
$P_{\a}$ and $R$.
The denominator of this rational function is precisely the minimal multiindex coefficient
in the expansion of $P_{\a_0}(\rho^i_{z^\b w^\gamma}(0,\bar\chi,0))$,
and hence it does not vanish since we have assumed this coefficient to be nonzero.

The final observation involves sufficiently large truncations of
the series in $\6R^t_0$ (similar to \cite{Fo04}).
It is clear that the dimension of the corresponding truncation space for the coefficients 
of $P_{\a_0}(\rho^i_{z^\b w^\gamma}(0,\bar\chi,0))$ is arbitrarily large,
whereas the dimension of the polynomial coefficients of $P_{\a}$ and $R$ is fixed by our choice.
Hence, choosing sufficiently large truncations,
the condition for the coefficients of $\rho^{i_0}_{z^{\a_0}}(0,\bar\chi,0)$
to be given by a rational function as above, defines a nowhere dense subset.
Going back to the space $\6R^t_0$ before the truncation,
we can see that the corresponding subset there is also nowhere dense,
hence is of the first category as desired. The details are left to the reader.
\epf

\appendix

\section{Obstructions to biholomorphic equivalence to real-algebraic submanifolds}\Label{sec-eq}

Here we briefly illustrate how our methods can be used to obtain obstructions to biholomorphic equivalence to real-algebraic submanifolds. The proofs are self-contained and do not depend on the previous sections. In particular, no elaborate weight estimates are needed here.

\bt\Label{Q-equiv}
Let $M\subset \C^{n+d}$ be a real-analytic generic submanifold of codimension $d$ through $0$
given by an equation $w=Q(z,\bar z,\bar w)$, $(z,w)\in\C^n\times\C^d$, where $Q$ is a $\C^d$-valued holomorphic function in a neighborhood of the origin.
Suppose that  $M$ is biholomorphically equivalent to a real-algebraic generic submanifold of $\C^{n+d}$.
Then any set of $n+1$ partial derivatives of the components,
$Q^{i_1}_{z^{\a_1}}(0,\bar\zeta),\ldots, Q^{i_{n+1}}_{z^{\a_{n+1}}}(0,\bar\zeta)$, with $\zeta$ varying in the Segre varity $S_0$, is algebraically dependent, i.e.\ satisfies a nontrivial polynomial equation
$P(Q^{i_1}_{z^{\a_1}}(0,\bar\zeta),\ldots, Q^{i_{n+1}}_{z^{\a_{n+1}}}(0,\bar\zeta))=0$.
\et

\bpf
Let $H=(F,G)\colon (\C^{n}\times\C^d,0)\to(\C^{n}\times\C^d,0)$ be a local biholomorphic map sending a neighborhood of $0$ in $M$ into a real-algebraic generic submanifold $M'\subset \C^{n+d}$ that we may assume being given by 
$w'=Q'(z',\bar z',\bar w')$, where $Q'$ is a (complex{-})algebraic  holomorphic $\C^d$-valued function satisfying $Q'_{z'}(0,0,0)=0$. The latter implies that the $d\times d$ matrix $G_w(0)$ is invertible. Then for every $k\ge1$, the holomorphic map 
\begin{equation}\Label{mu}
\mu_{Q',k}\colon\bar\zeta'\in \1{S'_0}\mapsto (Q'^j_{z'^{\b}}(0,\bar\zeta'))_{1\le j\le d, |\b|\le k}\in\C^N,
\end{equation}
with appropriate $N$,
is algebraic. Recall that $\dim S'_0=n$. Then by Chevalley's theorem (see e.g.\ \cite{M}, p.~72), the image of $\mu_{Q',k}$ is contained in an algebraic variety of dimension $n$.

The property that $H$ sends $M$ into $M'$ can be expressed (after complexification) as
\begin{equation}\Label{send}
G(z,Q(z,\bar\zeta))=Q'(F(z,Q(z,\bar\zeta)),\bar H(\bar\zeta)).
\end{equation}
Differentiating in $z$ at $(z,\bar\zeta)\in \{0\}\times\1{S_0}$ and using the properties 
$Q(0,\bar\zeta)=0$  for $\zeta\in S_0$, we conclude by induction on $|\a|$ that each derivative $Q^i_{z^{\a}}(0,\bar\zeta)$, $1\le i\le d$, $|\a|\le k$,
can be expressed as a rational function of $Q'^j_{z'^{\b}}(0,\bar H(\bar\zeta))$, $1\le j\le d$, $|\b|\le k$, with poles away from $\mu_{Q',k}(0)$. In particular, applying Chevalley's theorem we see that, for 
$\nu(\bar\zeta):=(Q^{i_1}_{z^{\a_1}}(0,\bar\zeta),\ldots, Q^{i_{n+1}}_{z^{\a_{n+1}}}(0,\bar\zeta))$
 the image of $\nu\circ \bar H$ is also contained in in an algebraic variety of dimension $n$.
Furthermore, since $\bar H$ maps $\1{S_0}$ locally biholomorphically onto $\1{S'_0}$ near $0$, it follows that the image $\nu(\1{S_0})$ is contained in the same algebraic variety of dimension $n$. The claimed algebraic dependence now immediately  follows.
\epf

As in case of Theorem~\ref{main} above, we obtain an immediate consequence in the special case when $M$ is {\em rigid}:

\bc
Let $M\subset \C^{n+d}$ be a real-analytic generic submanifold of codimension $d$ through $0$
given in its rigid form by $\Im w=\phi(z,\bar z)$, $(z,w)\in\C^n\times\C^d$.
Suppose that  $M$ is biholomorphically equivalent to a real-algebraic generic submanifold of $\C^{n+d}$.
Then any set of $n+1$ partial derivatives of the components,
$\phi^{i_1}_{z^{\a_1}}(0,\bar\chi),\ldots, \phi^{i_{n+1}}_{z^{\a_{n+1}}}(0,\bar\chi)$ with $\chi\in\C^n$ near the origin, is algebraically dependent.
\ec

In particular, specializing further to the tube case, we have:

\bc\Label{tube}
Let $M\subset \C^{n+d}$ be a real-analytic generic submanifold of codimension $d$ through $0$
given in its tube form by $\Im w=\phi(\Im z)$, $(z,w)\in\C^n\times\C^d$.
Suppose that  $M$ is biholomorphically equivalent to a real-algebraic generic submanifold of $\C^{n+d}$.
Then any set of $n+1$ partial derivatives of the components,
$\phi^{i_1}_{z^{\a_1}}(x),\ldots, \phi^{i_{n+1}}_{z^{\a_{n+1}}}(x)$ with $x\in\R^n$ near the origin, is algebraically dependent.
\ec

Some of the algebraic dependence relations in Corollary~\ref{tube} (with $i_1=\ldots=i_{n+1}$, $|\a_s|\le 2$, in case $M$ is minimal and finitely nondegenerate and its infinitesimal CR automorphism algebra has minimum possible dimension) are contained in \cite{GM} as mentioned before in the introduction.

We conclude by mentioning that Proposition~\ref{def-segre} can be used to obtain a version of Theorem~\ref{Q-equiv} with algebraic dependence relations for the derivatives of (the components of) any defining function rather than the function $Q$, similarly to Theorem~\ref{main} (or Theorem~\ref{main1}) being obtained from Proposition~\ref{Q-relation}.

\end{document}